\documentclass[12pt,leqno]{article}

\usepackage[francais]{babel}
\usepackage{amsmath}
\usepackage{amssymb}

\usepackage[all]{xy}
\newtheorem{lemme}{Lemme}
\newtheorem{theo}[lemme]{Th\'eor\`eme}
\newtheorem{prop}[lemme]{Proposition}
\newtheorem{cor}[lemme]{Corollaire}
\newtheorem{defi}[lemme]{D\'efinition}

\newenvironment{rem}{\noindent{\sc
Remarque }\stepcounter{lemme}\thelemme .
\kern.2em}{\goodbreak\vskip10pt}

\def\demo{\medskip\goodbreak\noindent
     \hbox{\sc D\'emonstration : \kern .3em}\ignorespaces}%
  \def \qedbox{$\square$}%
  \def \qed{\ifmmode\qedbox
     \else\unskip\ \hglue0mm\hfill\qedbox\medskip
      \goodbreak\fi}%

\def\preuve{\medskip\goodbreak\noindent
     \hbox{\sc Preuve : \kern .3em}\ignorespaces}%
  \def \qedbox{$\square$}%
  \def \qed{\ifmmode\qedbox
     \else\unskip\ \hglue0mm\hfill\qedbox\medskip
      \goodbreak\fi}%
\def\endpreuve{\qed\goodbreak\vskip10pt}%

\def\Proj{\mathrm{Proj}}

\def\Hom{\mathrm{Hom}}

\def\id{\mathrm{id}}

\def\Mat{\mathrm{Mat}}
\def\mat{\mathrm{mat}}

\def\End{\mathrm{End}}

\def\tr{\mathrm{tr}}
\def\rg{\mathrm{rg}}
\def\mod{\mathrm{mod}}
\begin{document}

\centerline{\huge Quelques classes caract\'eristiques}
\vskip3truemm
\centerline{\huge en th\'eorie des nombres}

\vskip13mm

\centerline{Max Karoubi\footnote{ Universit\'e Denis
Diderot (Paris VII), UFR de Math\'ematiques, UMR 7586 du CNRS, case 7012,
2 place Jussieu, 75251 Paris Cedex 05, France (courriel :
karoubi@math.jussieu.fr).}  et Thierry
Lambre\footnote{Universit\'e Paris-Sud (Paris XI),
D\'epartement de Math\'ematiques, UMR 7586 et 8628 du CNRS, 91405
Orsay Cedex, France (courriel : thierry.lambre@math.u-psud.fr).}}

\vskip12mm



\noindent {\bf {Abstract :}}
Let $A$ be an arbitrary ring. We introduce
a Dennis trace map mod $n$, from $K_1(A;{\bf Z}/n)$ to the Hochschild homology
group with coefficients $HH_1(A;{\bf Z}/n)$. If $A$ is the ring of integers
in a number field, explicit elements of
$K_1(A,{\bf Z}/n)$ are constructed and the values of their Dennis trace mod
$n$ are computed. If $F$ is a quadratic field, we obtain this way non trivial
elements of the ideal class group of $A$. If $F$ is a cyclotomic 
field, this trace is
closely related to Kummer logarithmic derivatives; this trace leads to an
unexpected relationship between the first case of Fermat last theorem,
$K$-theory and the number of roots of Mirimanoff polynomials.
\bigskip

\noindent {\bf {R\'esum\'e : }}
Pour un anneau $A$ arbitraire, nous
cons\-trui\-sons  une trace de Dennis \`a
coefficients, de source le groupe de $K$-th\'eorie $K_1(A;{\bf Z}/n)$, de but
le groupe d'homologie
  de Hochschild $HH_1(A;{\bf Z}/n)$. Lorsque $A$ est l'anneau des entiers d'un
corps de nombres, nous construisons des \'el\'ements explicites  du groupe
$K_1(A;{\bf Z}/n)$ et
\'evaluons leurs traces de
Dennis. Dans le cas des corps quadratiques, nous en d\'eduisons des 
\'el\'ements
non triviaux du
groupe des classes de $A$. Dans le cas cyclotomique, cette trace est intimement
reli\'ee aux d\'eriv\'ees logarithmiques de Kummer, ce qui  permet 
de formuler un
lien inattendu entre le premier cas du dernier th\'eor\`eme de Fermat, la
$K$-th\'eorie et  le nombre de racines des polyn\^omes de
Mirimanoff.

\bigskip

\noindent {\bf {Mots-cl\'es :}}
K-th\'eorie \`a coefficients, trace de Dennis \`a coefficients, 
groupe des classes, corps
quadratiques et cyclotomiques, d\'eriv\'ee logarithmique.
\bigskip

\noindent {\bf {Classifications AMS 1991 :}}
Primaire 11 R 29, 19 D 55.
  Secondaire 11 R 18, 18 F 30, 19 F 99

\bigskip

\noindent 1. La trace de Dennis \`a coefficients.

\noindent 1.1. Le groupe $K_1(A;{\bf Z}/n)$.

\noindent 1.2. L'alg\`ebre diff\'erentielle gradu\'ee $\Omega^*(A)$.

\noindent 1.3. Le module gradu\'e $HH_*(A;{\bf Z}/n)$.

\noindent 1.4. Calcul diff\'erentiel non commutatif d'ordre un.

\noindent 1.5. La trace de Dennis $D_1^{(n)}$.

\noindent 1.6. Le cas des alg\`ebres commutatives.

\noindent 1.7. Les traces d'ordre sup\'erieur.

\medskip
\noindent 2. \'Etude de l'anneau des entiers d'un corps de nombres.

\noindent 2.1. Description de $K_1(A;{\bf Z}/n)$ en termes d'id\'eaux.

\noindent 2.2. Le lemme $(N,N_1)$ de construction d'\'el\'ements de
$K_1(A;{\bf Z}/n)$.

\noindent 2.3. Description de $K_1(A;{\bf Z}/n)$ en termes d'ad\`eles.

\noindent 2.4. Description de $\Omega_{dR}^1(A)/(n)$.

\noindent 2.5. Description de la trace de Dennis \`a coefficients.

\medskip

\noindent 3. Applications aux corps de nombres de petit degr\'e.

\noindent 3.1. Un th\'eor\`eme de Y. Yamamoto.

\noindent 3.2. Construction d'\'el\'ements non triviaux du groupe des
classes.

\noindent 3.3. Exemples de $n$-torsion du groupe des classes : cas d'un corps
quadratique imaginaire.

\noindent 3.4. Exemples de $n$-torsion ramifi\'ee du groupe des
classes: cas d'un corps quadratique.
  \medskip

\noindent 4. Applications \`a la cyclotomie.

\noindent 4.1. Notations et strat\'egie g\'en\'erale.

\noindent 4.2. Emploi du groupe $K_1(A/p;{\bf Z}/p)$.

\noindent 4.3. Emploi du groupe $K_1(R;{\bf Z}/p)$.

\noindent 4.4. Lien avec les d\'eriv\'ees logarithmiques de Kummer.

\noindent 4.5. Application au premier cas du dernier th\'eor\`eme de Fermat.

\noindent 4.6. Lien avec les nombres de Bernoulli.

\bigskip
Dans tout ce texte, nous utiliserons les conventions suivantes.
Soit $M$ un groupe et $n$ un entier naturel.
On d\'efinit $\cdot n:M\to M$ par
$\cdot n(z)=nz$ si le groupe est not\'e additivement
et par $\cdot n(z)=z^n$ si le groupe est not\'e
multiplicativement.
On pose $M_{(n)}=\ker(\cdot n)$,
$M^{(n)}=\hbox {Im}(\cdot n)$
et $M/(n)=M/M^{(n)}$.\break
Tous les anneaux sont suppos\'es
unitaires.
Le groupe des  unit\'es d'un anneau $A$ est not\'e $A^\times$. La 
classe mod. $n$
d'un entier $a\in {\bf Z}$ est not\'ee $\overline a$.

\bigskip
\noindent {\bf \large Introduction : }

Soit $k$ un anneau commutatif unitaire et soit $A$ une $k$-alg\`ebre
(associative, unitaire). K. Dennis [10], A. Connes [7] et M. Karoubi [11] ont
cons\-truit des homomorphismes ``classes caract\'eristiques''
$D_i:K_i(A)\to HH_i(A)$, ($i\geq 0$) et $ch_{i,\ell}:K_i(A)\to
HC_{i+2\ell}(A)$ , ($i\geq 0$ et $\ell\geq 0$), de source la $K$-th\'eorie
de $A$, de but l'homologie de Hochschild de $A$ (pour la trace de Dennis
$D_i$) ou  l'homologie cyclique  (pour les caract\`eres de Chern
$ch_{i,\ell}$).

Lorsque $A$ est l'anneau des entiers d'un corps de nombres, ces classes
carac\-t\'e\-ristiques sont inop\'erantes, car trop souvent triviales (d\`es
que $i$ est pair, [14]). En rempla\c cant la
$K$-th\'eorie usuelle par la $K$-th\'eorie \`a coefficients ${\bf 
Z}/n$, nous construisons pour
$i=1$ une trace de Dennis \`a coefficients $\overline D_1:K_1(A;{\bf 
Z}/n)\to HH_1(A;{\bf
Z}/n)$. Pour un anneau de Dedekind $A$, nous explicitons le groupe $K_1(A;{\bf
Z}/n)$ ainsi que cette trace $\overline D_1$. Celle-ci s'av\`ere 
suffisamment riche pour
d\'etecter sous certaines conditions  des \'el\'ements de $n$-torsion du groupe
des classes.
  Des exemples d\'etaill\'es sont propos\'es
pour les corps quadratiques. Dans le cas d'un corps cyclotomique, ces classes
carac\-t\'e\-ristiques se trouvent reli\'ees de mani\`ere  inattendue aux
d\'eriv\'ees logarithmiques de Kummer et aux
polyn\^omes introduits au d\'ebut du si\`ecle par M. Mirimanoff lors de ses
recher\-ches sur le dernier th\'eor\`eme de Fermat.

   Les r\'esultats principaux de cet article sont les suivants.

\begin{theo}
Soit $A$ un anneau unitaire et soit $n\geq 2$ un entier naturel. Il 
existe un morphisme
naturel de groupes $\overline D_1$ s'ins\'erant dans le diagramme
commutatif suivant o\`u les lignes sont exactes.

\newsavebox{\maboiteun}
\savebox{\maboiteun}(350,80)[]{\scriptsize
$
\xymatrix{
K_1(A)\ar[d]^-{D_1} \ar[r]^-{.n} &K_1(A) \ar[d]^-{D_1} \ar[r]^-\rho &
K_1(A;{\bf Z}/n) \ar[d]^-{\overline D_1} \ar[r]^-{\partial} & K_0(A)
\ar[d]^-{D_0}\ar[r]^-{.n} &
K_0(A)\ar[d]^-{D_0} \ar[r] &K_0(A)/(n)\ar[d]^-{\overline D_0}\\
HH_1(A)\ar[r]^-{.n} & HH_1(A)\ar[r]^-{\rho} & HH_1(A;{\bf Z}/n) 
\ar[r]^-{\partial}
& HH_0(A) \ar[r]^-{.n} & HH_0(A)\ar[r] & HH_0(A;{\bf Z}/n)
  }
$}
\usebox{\maboiteun}
\end{theo}

Soit $A$ un anneau de Dedekind  de corps des fractions $F$. On d\'esigne par
$I(A)$le mono\"{\i}de des id\'eaux fractionnaires de $A$. Posons
$${\cal U}(A;{\bf Z}/n):=\{x\in F^\times/(n)\mid\exists\ I\in I(A),\
xA=I^n\}\cdot$$

\begin{theo}
Le groupe $K_1(A;{\bf Z}/n)$ d'un anneau de Dedekind $A$ est 
isomorphe au groupe
${\cal U}(A;{\bf Z}/n)\oplus SK_1(A)/(n)\cdot$
\end{theo}

Lorsque  $A$ est l'anneau des entiers d'un corps de nombres $F$,
le lemme suivant permet de construire des \'el\'ements du groupe
$K_1(A;{\bf Z}/n)$.

\begin{lemme} (``lemme $N$-$N_1$'') Soit $A$ l'anneau des entiers d' 
un corps de
nombres
$F$. Notons
$Cl(A)$ le groupe des classes de $A$.  Soit $u$ un \'el\'ement non nul  de $A$.
Pour que l'\'el\'ement
$[u]$ de
$F^\times/(n)$ appartienne \`a $K_1(A;{\bf Z}/n)$, il suffit que la 
norme $N(u)$
soit une puissance $n$-i\`eme dans
${\bf Z}$ et que
$(N(u), N_1(u))=1$, $N_1(u)$ \'etant  le coefficient de $X$ dans le
polyn\^ome caract\'eristique de $u$, consid\'er\'e comme endomorphisme du
${\bf Q}$-espace vectoriel $F$.
\end{lemme}

En appliquant ces deux r\'esultats aux corps quadratiques, nous aboutissons au
\begin{theo}
Soit $A$ l'anneau des entiers d'un corps de nombres quadratique $F$. 
Soit $n$ un diviseur
impair du discriminant $\delta$ de $F$. Si $F$ est r\'eel, on suppose 
que l'unit\'e
fondamentale
$\varepsilon =\displaystyle\frac{\varepsilon_1+\varepsilon_2\sqrt{\delta}}{2}$
est telle que
$n$ divise $\varepsilon_2$.
Alors, il existe une ``classe caract\'eristique secondaire''
$$d_1^{(n)}:Cl(A)_{(n)}\to{\bf Z}/n,$$ non triviale en g\'en\'eral.
\end{theo}

Dans le cas o\`u $A={\bf Z}[\zeta]$ est l'anneau des entiers du corps 
cyclotomique
${\bf Q}[\zeta_p] $ (avec $p$ impair), la conjugaison complexe scinde le groupe
$K_1(A;{\bf Z}/p)$ en deux sous-groupes dont la partie antisym\'etrique est
not\'ee
$K_1^-(A;{\bf Z}/p)$
et s'ins\'ere dans la suite exacte courte
$$\xymatrix{1\ar[r]&\mu_p\ar[r]&K_1^-(A;{\bf Z}/p)\ar[r]&Cl(A)_{(p)}^-\ar[r]&1
}$$ o\`u $Cl(A)$ est le groupe des classes de $A$.
Posons $$d_p^-=\dim_{{\bf Z}/p}Cl(A)_{(p)}^-=
\dim_{{\bf Z}/p}K_1^-(A;{\bf Z}/p)-1\cdot$$ Les
nombres premiers r\'eguliers sont caract\'eris\'es par
$d_p^-=0$. La th\'eorie analytique des nombres montre que pour $p\geq 7$, on
a
$d_p^-\leq (p-1)/4\cdot$

Si $(p,a,b,c)$ satisfont aux hypoth\`eses du premier cas du dernier 
th\'eor\`eme de
Fermat (c'est-\`a-dire si $p$  premier impair,  $a$, $b$, $c$ trois 
entiers tels que
$a^p=b^p+c^p$ avec $(a,b,c)=1$ et $p$ ne divisant pas $abc$),
l'\'el\'ement $$z=\frac{a-b\zeta}{a-b\zeta^{-1}}\ \mod\  \ F^{\times(p)}$$ de
$F^\times/(p)$ appartient \`a $K_1^-(A;{\bf Z}/p)$. En calculant la trace de
Dennis
\`a coefficients de cet \'el\'ement on aboutit \`a la minoration $d_p\geq 1$,
ce qui prouve le th\'eor\`eme de Kummer sur les nombres premiers r\'eguliers.
Les coefficients de cette trace de Dennis sont d'ailleurs li\'es aux 
d\'eriv\'ees
logarithmiques de Kummer. Cette trace permet
\'egalement d'obtenir une autre minoration de
$d_p$, qui s'exprime \`a l'aide des polyn\^omes de Mirimanoff
$M_k(X)\in{\bf Z}/p[X]$, ($1\leq k\leq p-1$), d\'efinis par
$M_k(X)=\sum_{j=1}^{p-1}j^{k-1}X^j$. Pour $t\in {\bf Z}/p\setminus\{0,1,-1\}$,
posons
$$r_p(t):=\sharp\{k,\ 1 \leq k\leq \frac{p-1}{2},\ M_{2k+1}(t)\not=0\}$$
et $$r_p=\min\{ r_p(t),\ t\in {\bf Z}/p\setminus\{ 0,1,-1\}\}\cdot$$
Les nombres $r_p(t)$ et $r_p$   sont li\'es aux
congruences de Kummer et \`a la divisibilit\'e des nombres de Bernoulli.
Nous montrons la minoration suivante de $d_p^-$, tr\`es voisine d'une 
minoration
obtenue pr\'ecedemment par Br\" uckner \`a l'aide  de d\'eveloppements de
fonctions logarithmes.
\begin{theo}
Si $(p,a,b,c)$ satisfont aux hypoth\`eses du premier cas du dernier
th\'eor\`eme de
Fermat, alors on a
$$d_p^-\geq r_p(a/b)-2.$$
\end{theo}

Nous en d\'eduisons que les nombres $r_p(t)$ et $r_p$  sont 
\'egalement li\'es au
th\'eor\`eme de Fermat-Wiles sous la forme suivante.
\begin{cor}
Soit $p$ un nombre premier. Si  $r_p\geq (p+11)/4$, alors le premier
cas du dernier th\'eor\`eme de Fermat est satisfait pour $p$.
\end{cor}

\medskip
\noindent {\sl Remerciements}:
  Karim Belabas (Universit\'e Paris-Sud)  et Thong Nguyen Quang Do (Universit\'e
de Franche-Comt\'e) ont  bien volontiers accept\'e de lire une version
pr\'eliminaire de certains passages de ce texte.  Nous sommes heureux de les
remercier pour leur lecture attentive. Nous remercions \'egalement le
rapporteur anonyme de la Note aux C.R. Acd.Sci.  ([12]) pour nous 
avoir signal\'e la
r\'ef\'erence [6].
Signalons enfin que dans le cas des anneaux commutatifs,  J. Berrick (National
University of Singapore), a construit ([4]) des invariants voisins des notres
au moyen de ``matrices entrelac\'ees''.

\section{La trace de Dennis \`a coefficients.}

Apr\`es avoir rappel\'e en 1.1 les propri\'et\'es du groupe 
$K_1(A;{\bf Z}/n)$, on d\'efinit en 1.3
le groupe
$HH_1(A;{\bf Z}/n)$. La trace de Dennis \`a coefficients n\'ecessite 
une notion adapt\'ee de trace
pour une $1$-forme diff\'erentielle non commutative. Cette 
construction est d\'etaill\'ee en 1.4 et
la trace de Dennis \`a coefficients est construite en 1.5.
Si l'anneau consid\'er\'e est commutatif, une construction plus 
\'el\'ementaire est propos\'ee en
1.6 gr\^ace aux formes diff\'erentielles \`a la de Rham. Des traces 
d'ordres sup\'erieurs,
de source
$K_i(A;{\bf Z}/n)$, de but $HC_i^-(A;{\bf Z}/n)$  ou
$HH_i(A;{\bf Z}/n)$ ($i\geq 1$), sont d\'efinies en 1.7. Ces 
d\'efinitions 1.7 ne sont
utilis\'ees nulle part ailleurs dans le texte.

\subsection{Le groupe $K_1(A;{\bf Z}/n)$.}

Soit $A$ un anneau. On note
$Proj(A)$ la cat\'egorie des $A$-modules \`a droite, projectifs et de 
type fini.
Soit $P\in Proj(A)$ et
$\alpha\in Aut_A(P)$. Rappelons que le groupe de Bass
$K_1(A)$ peut-\^etre vu comme le quotient du groupe ab\'elien libre 
engendr\'e par les
classes d'isomorphie de paires
$(P,\alpha)$ modulo le sous-groupe engendr\'e par les deux types 
d'\'el\'ements suivants:

(a) $(P,\alpha)+(P,\beta)-(P,\alpha\beta)$

(b) $(P_1\oplus P_2,\alpha_1\oplus\alpha_2)-(P_1,\alpha_1)-(P_2,\alpha_2).$

On note $[P,\alpha]\in K_1(A)$ la classe de $(P,\alpha)$

Soient $A$ et $B$ deux anneaux unitaires et soit $\varphi:Proj(A)\to 
Proj(B)$ un
foncteur additif. Le foncteur $\varphi$ est dit cofinal ([2], VII.1, p. 345)
si tout objet $R$ de $Proj(B)$ est facteur direct d'un objet de la forme
$\varphi (P)$ avec $P$ objet de $Proj(A)$, c'est-\`a-dire 
$\varphi(P)\cong R \oplus S$ avec $S\in
Ob(Proj(B))$.

\`A un tel foncteur cofinal $\varphi$, on associe la cat\'egorie
${\cal C}(\varphi)$ dont les objets sont les
triplets $(P,\alpha ,Q)$ avec $P\in Ob(Proj(A))$,
$Q\in Ob(Proj(A))$, et o\`u $\alpha$ est un isomorphisme de
$B$-modules $\varphi(P)\cong\varphi(Q)$.

Un morphisme dans ${\cal C}(\varphi)$ de source
$(P,\alpha,Q)$, de but
$(P_1,\alpha_1,Q_1)$ est un couple $(f,g)$ de morphismes de
$A$-modules $f:P\rightarrow P_1$ et
$g:Q\rightarrow Q_1$
tels que
$\varphi(g)\circ\alpha=\alpha_1\circ\varphi(f)$.
L'ensemble des classes d'isomorphie d'objets de la cat\'egorie
${\cal C}(\varphi)$ est un mono\"{\i}de ab\'elien; on note
$K({\cal C}(\varphi))$ le groupe de Grothendieck associ\'e.

\begin{defi}
Le quotient de $K({\cal C}(\varphi))$ par le sous-groupe $N$ engendr\'e par les
\'el\'ements
$$(P,\alpha, Q)+(Q,\beta, R)-(P,\beta\alpha, R)$$
est not\'e $K(\varphi)$.
On d\'esigne par $[P,\alpha,Q]$ la classe de
$(P,\alpha,Q)\in Ob({\cal C}).$

\end{defi}

\begin{rem}
Soit $\alpha':P\rightarrow Q$ un isomorphisme de $A$-modules. Alors
l'\'el\'ement $x=[P,\varphi(\alpha'),Q]\in K(\varphi)$ est nul.
En effet $(\alpha',\id_Q)$ est un isomorphisme de la cat\'egorie ${\cal C}$ de source 

$[P,\varphi(\alpha'),Q]$, de but $[Q,\id_{\varphi(Q)},Q]$, donc
$x=[Q,\id_{\varphi(Q)},Q]=0$. On en d\'eduit qu'on peut toujours supposer qu'un
\'el\'ement de $K(\varphi)$ est de la forme $[P,\alpha, L]$ avec $L$ 
libre car si
$[P,\alpha, Q]$ appartient \`a $K(\varphi)$ et si $Q\oplus R$ est 
isomorphe \`a un module
libre $L$, on a $[P,\alpha, Q]=[P\oplus R,\alpha\oplus\id_{nR},L]$.
\end{rem}

D'apr\`es [2], VII.5, p. 375, on a
\begin{theo}
Soient $A$ et $B$ deux anneaux unitaires et soit
$\varphi:Proj(A)\rightarrow Proj(B)$ un foncteur cofinal.
Alors, on a
une suite exacte de groupes ab\'eliens
$$
\xymatrix{
K_1(A) \ar[r]^-{\varphi_1} &K_1(B) \ar[r]^-{\rho} &K(\varphi) 
\ar[r]^-{\partial}
&K_0(A)\ar[r]^-{\varphi_0} &K_0(B)
}
$$
o\`u les applications $\varphi_1$, $\partial$ et $\varphi_0$
sont d\'efinies ainsi
$$\varphi_1([P,\alpha])=[\varphi(P),\varphi(\alpha)],$$
$$\partial([P,\alpha,Q])=[P]-[Q],$$
$$\varphi_0([P])=[\varphi(P)].$$
Pour $[R,\lambda]\in K_1(B)$, l'application $\rho$ est d\'efinie  par 
cofinalit\'e
de $\varphi$ en \'ecrivant $R\oplus S\cong\varphi(P)$ et en posant
$\rho([R,\lambda])=[P,\lambda\oplus id_{S},P].$

\end{theo}

\begin{rem} De cette suite exacte, on d\'eduit que le groupe 
$K(\varphi)$ est l'extension
$$
\xymatrix{
1 \ar[r]  &K_1(A)/ \varphi_1(K_1(A)) \ar[r]^-{\rho} & K(\varphi)
\ar[r]^-{\partial}  & Ker(\varphi_0)\ar[r] & 0\, .
}
$$
\end{rem}
  Appliquons cette construction au contexte suivant.
  \begin{defi}
Soit $n$ un entier, $n\geq 2$. Le foncteur
  $.n:Proj(A)\rightarrow Proj(A)$ est d\'efini par
  $.n(P)=nP=P\oplus\cdots \oplus P$
  et
  $.n(f)=nf=f\oplus \cdots \oplus f$
  ($n$ facteurs).
  Le groupe $K(.n)$ est not\'e $K_1(A;{\bf Z}/n)$.
  Un \'el\'ement de $K_1(A;{\bf Z}/n)$ est de la forme
  $[P,\alpha,Q]$, o\`u
  $P$ et $Q$ sont dans $Proj (A)$ et o\`u $\alpha$ est un isomorphisme 
de $A$-modules $nP\cong nQ$.
\end{defi}
On a par cons\'equent la suite exacte de $K$-th\'eorie \`a coefficients
  $$
\xymatrix{
K_1(A) \ar[r]^-{.n} &K_1(A) \ar[r]^-{\rho} &K_1(A;{\bf Z}/n) \ar[r]^-{\partial}
&K_0(A)\ar[r]^-{.n} &K_0(A)
}
$$

  L'extension de la remarque
  ci-dessus s'\'ecrit
  $$
\xymatrix{
1 \ar[r]  &K_1(A)/ (n) \ar[r]^-{\rho} & K_1(A;{\bf Z}/n)
\ar[r]^-{\partial}  & K_0(A)_{(n)}\ar[r] & 0\,
}\leqno{(\ddagger )}
$$
avec $\rho([P,\alpha])=[P,\alpha\oplus id_{(n-1)P},P]$
et
$\partial ([P,\alpha,Q])=[P]-[Q].$
Cette extension fournit le premier exemple de calcul de $K$-th\'eorie 
\`a coefficients.
\medskip

\noindent
{\bf Exemple:} Soit $A$ un anneau commutatif local. Pour $n\geq 2$, le groupe
$K_1(A;{\bf Z}/n)$ est \'egal \`a $A^\times/(n)$.

\medskip
Revenons au cas g\'en\'eral. On a les relations int\'eressantes suivantes.

\begin{lemme} Soit $A$ un anneau et
soit $n$ un entier, $n >1$.

\noindent Si $n \not\equiv 2\  \mod \ 4 $\ ,\ on a $nK_1(A;{\bf Z}/n) = 0$\, ;

\noindent et si $n \equiv 2 \ \mod \ 4$\ ,\  on a $2n K_1(A;{\bf Z}/n) = 0$\, .

\noindent En particulier, si $p$ est un nombre premier impair, 
$K_1(A,{\bf Z}/p)$
est un ${\bf Z}/p $-espace vectoriel.
\end{lemme}

\preuve montrons d'abord que pour tout $n>1$, on a $2n K_1(A; {\bf
Z}/n)=0$. Si $p$ et $q$ sont deux entiers et si $P \in  \Proj(A)$, on a un
isomorphisme $p(qP) \cong q(pP)$. Pour $q = p = n$, notons $v_P : n^2P
\cong n^2P$ cet isomorphisme. La relation \'evidente $v^2_P = \id$
montre que l'\'el\'ement $x_P = [P, v_P,P]$ de $K_1(A;{\bf Z}/n)$ est tel
que $2x_P = 0$. Soit $x = [P,\alpha,Q]\in K_1(A;{\bf Z}/n)$ ; on a $nx =
[nP,\beta,nQ]$ avec $\beta = v_Q\circ n\alpha \circ v_P$, ce qui donne $nx
= x_Q + x' + x_P$ avec $x' = [nP, n\alpha, nQ]=0$. On en d\'eduit $2nx = 0$,
ce qui montre $2n K_1(A;{\bf Z}/n)=0$.

Par ailleurs, la suite exacte de $K$-th\'eorie \`a coefficients donne
$n^2K_1(A;{\bf Z}/n)$ $=0$. Pour
$n$ impair, les relations $2nK_1(A;{\bf Z}/n)=0$ et $n^2K_1(A;{\bf
Z}/n)=0$ conduisent \`a $nK_1(A;{\bf Z}/n)=0$.

Il reste \`a montrer que pour $n \equiv 0 \ \mod \ 4$, on a \'egalement
$nK_1(A;{\bf Z}/n)=~0$. Fixons une notation : pour $\sigma \in
\mathfrak{S}_n$ et $P \in  \Proj(A)$, on note $\sigma$ l'automorphisme
de $nP$ d\'efini pour $(z_1,\ldots ,z_n)\in nP$ par $\sigma(z_1,\ldots
,z_n) = \break z_{\sigma(1)},\ldots ,z_{\sigma(n)}\in nP$. On v\'erifie
facilement que pour tout $n >0$ et tout $P \in  \Proj(A)$, l'isomorphisme
$v_P : n(nP) \cong n(nP)$ est le produit de $\displaystyle \frac{n(n-1)}2$
transpositions \`a supports disjoints. Ce nombre est pair si $n \equiv 0$
ou $1 \ \mod \ 4$. La proposition suivante montre que dans  ce cas,
l'\'el\'ement $[P,v_P,P]$ de\break $K_1(A;{\bf Z}/n)$ est nul et par suite que
$nK_1(A;{\bf Z}/n)=0$.
\endpreuve

\begin{prop}
Soit $n \geq 4$ et $\tau\in \mathfrak{S}_n$ le produit de deux
transpositions \`a supports disjoints. Alors pour tout $P \in  \Proj(A)$
l'\'el\'ement $[P,\tau,P]$ de $K_1(A; {\bf Z}/n)$ est nul.
\end{prop}

\preuve dans $\mathfrak{S}_n$, $\tau$ est conjugu\'e \`a $\tau' = (12)
(34)$. Dans $K_1(A; {\bf Z}/n)$, on a donc $[P,\tau,P] = [P,\tau',P]$. La
matrice $T' \in \mathrm{GL}({\bf Z})$ de $\tau'$ appartient \`a\hfill\break
$[\mathrm{GL\,}({\bf Z}),
\mathrm{GL\,}({\bf Z})]$
donc
$[P,\tau',P]=0$.
\endpreuve

\subsection{L'alg\`ebre diff\'erentielle gradu\'ee $\Omega^*_{nc}(A)$.}

Soit $k$ un anneau commutatif unitaire et soit $A$ une $k$-alg\`ebre
unitaire. D\'esi\-gnons par $\displaystyle \mu : A \otimes_k A \to  A$ la
multiplication de $A$. Posons $\Omega^1_{nc}(A) := \ker \mu$. On sait que
$\Omega^1_{nc}(A)$ est le sous-bimodule de $\displaystyle A \otimes_k
A$ engendr\'e par $$\{1 \otimes a - a \otimes 1, a \in  A\}\cdot$$ On introduit
la 1-forme diff\'erentielle non commutative $$d_{nc}a := 1\otimes a- 
a\otimes 1.$$ La
structure de
$A$-bimodule \'evident de $\Omega^1_{nc}(A)$ se lit en terme de formes
diff\'erentielles par la relation
$$
d_{nc}(a_1) \cdot a_2 = d_{nc}(a_1a_2) - a_1 d_{nc}a_2 \leqno(*)
$$
de sorte que tout \'el\'ement de $\Omega^1_{nc}(A)$ s'\'ecrit comme somme
de 1-formes diff\'eren\-tielles $a_0d_{nc}a_1$.
En tant que $k$-module, on a $\Omega^1_{nc}(A)\cong A\otimes_k(A/k)$.

A partir de $\Omega^1_{nc}(A)$, on construit une $k$-alg\`ebre
diff\'erentielle gradu\'ee $\Omega^*_{nc}(A)$ en posant
$\displaystyle \Omega^r_{nc}(A) :=\Omega^1_{nc}(A)\otimes_A \ldots
\otimes_A \Omega^1_{nc}(A)$ ($r$ facteurs). La structure de $A$-module
\`a droite sur $\Omega^r_{nc}(A)$ se d\'efinit en jouant \`a saute-moutons \`a
l'aide de la formule $(*)$ :
$$
a_0d_{nc}a_1\ldots d_{nc}a_r \cdot \alpha = a_0d_{nc}a_1\ldots d_{nc}(a_r
\alpha)-a_0d_{nc}a_1\ldots d_{nc}a_{r-1}\cdot a_r d_{nc}\alpha = 
\mathrm{etc}\,.
$$
Tout \'el\'ement de $\Omega^r_{nc}(A)$ est somme de $r$-formes
diff\'erentielles du type
$$a_0d_{nc}a_1\ldots d_{nc}a_r.$$ C'est encore la relation
$(*)$ qui d\'efinit le produit de $\omega_r \in  \Omega^r_{nc}(A)$ par
$\omega_s\in \Omega^s_{nc}(A)$.

En tant que $k$-module, on a $\Omega^r_{nc}(A) = A \otimes (A/k)^{\otimes
r}$ et le bord de Hochschild de $A$ admet une expression tr\`es simple
dans l'alg\`ebre $\Omega^*_{nc}(A)$~: pour $\omega_r\in \Omega^r_{nc}(A)$ avec
$\omega_r = \omega' d_{nc}\alpha\, , \, \alpha \in A$, on a $b(\omega_r) =
(-1)^{r-1}(\omega' \alpha - \alpha \omega')$. Le complexe de
Hochschild normalis\'e, traditionnellement not\'e $(\overline 
C_*(A),b)$ n'est autre que
$(\Omega^*_{nc}(A),b)$. Pour $\omega\in \Omega^*_{nc}(A)$ tel que 
$b(\omega)=0$,
on note $[\omega]$ sa classe d'homologie de Hochschild.

\subsection{Le $k$-module gradu\'e $HH_*(A;{\bf Z}/n)$.}

Soit $k$  un anneau commutatif unitaire, soient $(C_*,d)$ et $(C'_*, d')$
des complexes (de cha\^\i nes) de $k$-modules et soit $f : (C_*,d)\to
(C'_*, d')$ un morphisme de complexes. Le c\^one de $f$ est le complexe
de cha\^\i nes $(\mathrm{co}(f),\partial )$ d\'efini par $\mathrm{co}(f)_r
= C'_r
\oplus C_{r-1}$ et $\partial (a',a) = (d'(a')+f(a), -d(a))$. De la suite exacte
courte de complexes
$$
0 \to  (C'_*,d') \to  (\mathrm{co}(f),\partial ) \to  (C_*,d)_{[1]}\to 0
$$
on d\'eduit une suite exacte longue d'homologie dont le con\-nec\-tant\break
$H_{r}\left( (C_*,d)_{[1]}\right) \to  H_{r-1}(C'_*,d')$ n'est autre 
que $H_{r-1}(f)$.

\begin{defi}
Soient $k$ un anneau commutatif unitaire et $A$ une $k$-alg\`ebre
unitaire. Soit $(\Omega_*(A),b)$ le complexe de Hochschild de $A$.
L'homologie du c\^one de l'application $.n : (\Omega_*(A),b) \to 
(\Omega_*(A),b)$
d\'efinie par $.n(a_0\otimes a_1\otimes \ldots \otimes a_r) =
na_0\otimes a_1\otimes \ldots \otimes a_r$ s'appelle l'homologie de
Hochschild \`a coefficients ${\bf Z}/n$ de l'alg\`ebre $A$. On note cette
homologie $HH_*(A; {\bf Z}/n)$.
\end{defi}

Pour $(\omega_r, \omega_{r-1})\in \Omega^r(A)\oplus \Omega^{r-1}(A)$
telle que $\partial (\omega_r, \omega_{r-1})=0$, on note $[\omega_r,
\omega_{r-1}]$ sa classe d'homologie de Hochschild \`a coefficients.

De cette d\'efinition de $HH_*(A; {\bf Z}/n)$, on extrait la suite exacte
longue
$$\xymatrix{
\ldots \ar[r] &HH_1(A) \ar[r]^-{.n} & HH_1(A) \ar[r]^-{\rho} & HH_1(A;
{\bf Z}/n) \ar[r]^-{\partial} & HH_0(A) \\
& \ar[r]^-{.n} & HH_0(A) \ar[r] &  HH_0(A; {\bf Z}/n)\cdot
}
$$
avec $\rho([\omega_r]) = [\omega_r,0]$ et $\partial([\omega_r,
\omega_{r-1}]) = [\omega_{r-1}]$.

\subsection{Calcul diff\'erentiel non commutatif d'ordre 1.}

Soit $P$ un $A$-module projectif, \`a droite, de type fini. Un syst\`eme
de coordonn\'ees de $P$ (cf. [5], II.46) est une suite ${\cal S}=
(x_j, \varphi_j)_{1\leq j\leq r}$ avec $x_j\in P$, $\varphi_j \in  P^* =
\Hom_A(P,A)$ telle que pour tout $x \in  P$, on ait $ x =
\sum^r_{j=1}x_j\varphi_j(x)$.

Soit $u \in  \End P$. Par d\'efinition, la matrice de $u$ dans le syst\`eme
de coordonn\'ees ${\cal S}$ est la matrice $U = (U_{ij}) \in \Mat_{r,r}(A)$
d\'efinie par $U_{ij} =\varphi_i \circ u(x_j)$. La quantit\'e
$\sum^r_{i=1} \varphi_i \circ u(x_i)\in A$ s'appelle la
trace de $u$ relativement au syst\`eme de coordonn\'ees ${\cal S}$. On la
note $\tr(u, {\cal S})$, ou $\tr(u)$. On pose $\rg(P) = \rg(P,{\cal S}) =
\tr(\id_P, {\cal S}) = \sum^r_{i=1} \varphi_i(x_i)\in A$. Ces
quantit\'es d\'ependent  du syst\`eme de coordonn\'ees ${\cal S}$ choisi.

La connexion de Levi-Civita de $P$ (cf. [7] et [11]) est l'application $d_P
: P \to  P \otimes _A \Omega^1(A)$ dont l'expression dans le syst\`eme
de coordonn\'ees ${\cal S}$ est
$$
d_P \left( \sum^r_{j=1} x_j\varphi_j(x)\right) = \sum^r_{j=1} x_j
d_{nc}(\varphi_j(x))\,.
$$
Comme toute connexion, $d_P$ n'est pas une application $A$-lin\'eaire mais
satisfait \`a  la relation :
$$
d_P(xa) = d_P(x) a + xd_{nc}a \ , \ x\in P \ , \ a \in  A\,.
$$

\begin{defi}
Soit $\alpha : P \to Q$ une application $A$-lin\'eaire \`a droite entre les
modules projectifs de type fini $P$ et $Q$. Soient $d_P$ et $d_Q$ les
connexions de Levi-Civita de $P$ et $Q$. L'application $A$-lin\'eaire
\begin{eqnarray*}
&&d_{nc}\alpha : P \to  Q \otimes _A \Omega^1(A) \text{ est d\'efinie par
}\\
&& d_{nc}\alpha = d_Q \circ \alpha - (\alpha \otimes \id)\circ d_P\,.
\end{eqnarray*}
\end{defi}

L'abus d'\'ecriture fr\'equent qui consiste \`a poser $d_P = d_Q = d_{nc}$
conduit \`a la formule $d_{nc} (\alpha(x)) = d_{nc}\alpha(x) + 
\alpha(d_{nc}x)$, $x \in P$.

Soient ${\cal S} = (x_j,\varphi_j)_{1\leq j\leq r}$ et ${\cal S'}= (y_i,
\psi _i)_{1\leq i\leq s}$ des syst\`emes de coordonn\'ees respectifs de
$P$ et $Q$. Il est facile de donner une interpr\'etation matricielle de
l'application $d\alpha$. Si $M = (M_{ij}) \in  \Mat_{s,r}(A)$ est la
matrice de $\alpha$, exprim\'ee dans les syst\`emes de coordonn\'ees
${\cal S}$ et ${\cal S'}$, avec $M_{ij} = \psi _i \circ \alpha(x_j)$, un
calcul simple montre que pour $x = \sum^r_{j=1} x_j\varphi_j(x)$, on a
$d_{nc}\alpha(x) = \sum_{i,j} y_i d(M_{ij}) \varphi_j(x)$, ce qui permet de
dire que $d_{nc}\alpha$ a pour matrice $d_{nc}M = (d_{nc} M_{ij})\in
\Mat_{s,r}(\Omega^1(A))$ dans les syst\`emes ${\cal S}$ et ${\cal S'}$.

Supposons de plus que $\alpha : P\to Q$ soit un isomorphisme de
$A$-modules. Soit $N \in  \Mat_{s,r}(A)$ la matrice de $\alpha^{-1}$
exprim\'ee dans les syst\`emes de coordonn\'ees ${\cal S'}$ et ${\cal S}$.
On a $MN = \mat_{{\cal S'}}(\id_Q)$ et $NM =\mat_{{\cal S}}(\id_P)$.
L'application compos\'ee
$\xymatrix{P \ar[r]^-{d_{nc}\alpha} & Q \otimes _A
\Omega^1(A) \ar[r]^-{\alpha^{-1}\otimes \id} & P \otimes _A \Omega^1(A))}$,
\'evidemment not\'ee\break $\alpha^{-1}d_{nc}\alpha$ est $A$-lin\'eaire, de
matrice
$NdM\in
\Mat_{r,r}(\Omega^1(A))$. La trace de cette matrice, $\tr(NdM)\in
\Omega^1(A)$, s'appelle la trace de $\alpha^{-1}d_{nc}\alpha$ dans les
syst\`emes de coordonn\'ees ${\cal S}$ et ${\cal S'}$.

\begin{prop}
Soient $P$ et $Q$ deux $A$-modules \`a droite, projectifs et de type fini,
de syst\`emes de coordonn\'ees respectifs ${\cal S}$ et ${\cal S'}$.
Posons $p = \rg(P,{\cal S})$ et $q = \rg(Q, {\cal S'})$. Soit $\alpha : P\to
Q$ un isomorphisme de $A$-modules et soit $\tr(\alpha^{-1}d_{nc}\alpha)\in
\Omega^1(A)$ la trace de $\alpha^{-1}d_{nc}\alpha$, exprim\'ee dans les
syst\`emes de coordonn\'es ${\cal S}$ et ${\cal S'}$.

Alors le bord de Hochschild de $\tr(\alpha^{-1}d_{nc}\alpha)$ est donn\'e par
$b(\tr(\alpha^{-1}d_{nc}\alpha)) = q - p\in A$.
\end{prop}

\preuve
on a $\tr(\alpha^{-1}d_{nc}\alpha)= \tr(Nd_{nc}M)$ d'o\`u

\noindent
$b(\tr(\alpha^{-1}d_{nc}\alpha))=b(\tr N d_{nc} M) = - \tr NM + \tr MN =
-\tr(\id_P, {\cal S}) +\tr(\id_Q,{\cal S'})\,.$
\endpreuve

\rem Dans la proposition ci-dessus, si $P=Q$ et si ${\cal S}={\cal 
S'}$, alors pour
tout automorphismes $\alpha$ de $P$, $\tr(\alpha^{-1}d_{nc}\alpha,{\cal S})$
est un cycle de Hochschild. Dans ce cas on note
$[\tr(\alpha^{-1}d_{nc}\alpha]\in HH_1(A)$, sa classe d'homo\-lo\-gie de
Hochs\-child.

\medskip

  Nous utiliserons \'egalement le r\'esultat suivant, de preuve facile
laiss\'ee au lecteur.

\begin{prop}
Soient $P$, $Q$ et $R$ des $A$-modules \`a droite, projectifs et de type fini ;
soient $\alpha\in \Hom_A(P,Q)$, $\beta\in \Hom_A(Q,R)$. Alors on a la
relation
$$
d_{nc}(\beta\circ \alpha) = d_{nc} \beta \circ \alpha + \beta \circ 
d_{nc}\alpha \,.
$$
En particulier $\alpha^{-1}d_{nc}\alpha = -d_{nc} \alpha^{-1}\circ \alpha$.
\end{prop}

\subsection{La trace de Dennis \`a coefficients $\overline D_1$.}

Pour tout entier $ r \geq  0$, K. Dennis ([10]) a construit un
morphisme
$D_r : K_r(A) \to HH_r(A)$.
Rappelons que si $r = 0$ et $[P]\in K_0(A)$,
on a $D_0([P]) = \rg(P,{\cal S})\ \mod\ [A,A]$ et que si $r=1$ et
$[P,\alpha]\in K_1(A)$, on a $D_1([P,\alpha]) = [\tr(\alpha^{-1}
d_{nc}\alpha)]\in HH_1(A)$.

\begin{theo}
Soient $A$ un anneau unitaire, $n$ un entier, $n \geq 2$ et soit $x =
[P_1,\alpha,P_2]\in K_1(A; {\bf Z}/n)$. Apr\`es avoir fix\'e des
syst\`emes de coordonn\'ees ${\cal S}_1$ et ${\cal S}_2$ sur $P_1$ et
$P_2$, on pose $p_1 = \rg(P_1,{\cal S}_1), p_2 = \rg(P_2,{\cal S}_2)$ et
$\overline D_1(x) = [\tr(\alpha^{-1}d_{nc}\alpha), p_1-p_2]\in 
HH_1(A; {\bf Z}/n)$.

Alors l'application $\overline D_1$ est un morphisme de groupes s'ins\'erant
dans le dia\-gramme commutatif \`a lignes exactes ci-dessous

\newsavebox{\maboitedeux}
\savebox{\maboitedeux}(350,60)[]{\scriptsize
$
\xymatrix{
K_1(A)\ar[d]^-{D_1} \ar[r]^-{.n} &K_1(A) \ar[d]^-{D_1} \ar[r]^-\rho &
K_1(A;{\bf Z}/n) \ar[d]^-{\overline D_1} \ar[r]^-{\partial} & K_0(A)
\ar[d]^-{D_0}\ar[r]^-{.n} &
K_0(A)\ar[d]^-{D_0} \ar[r] &K_0(A)/(n)\ar[d]^-{\overline D_0}\\
HH_1(A)\ar[r]^-{.n} & HH_1(A)\ar[r]^-{\rho} & HH_1(A;{\bf Z}/n) 
\ar[r]^-{\partial}
& HH_0(A) \ar[r]^-{.n} & HH_0(A)\ar[r] & HH_0(A;{\bf Z}/n)
  }
$}
\usebox{\maboitedeux}
\end{theo}

\preuve
la quantit\'e $c = \big(\tr(\alpha^{-1}d_{nc}\alpha), p_1-p-2\big)$ est
un cycle du c\^one de la multiplication
$.n: (\Omega^*(A),b)\to (\Omega^*(A),b)$. En effet,

$$\partial (c) =
b(\tr(\alpha^{-1}d_{nc}\alpha)) + n(p_1-p_2) ; $$
or
$$b(\tr(\alpha^{-1}d_{nc}\alpha)) = \rg(nP_2, n {\cal S}_2)-\rg(nP_1, n{\cal
S}_1),$$
  donc
$\partial (c)=0\,.$

Soit $[c] = [\tr \alpha^{-1} d_{nc}\alpha, p_1-p_2] \in  HH_1(A, {\bf Z}/n)$ la
classe d'homologie du cycle $c$. Cette classe est ind\'ependante du choix
du repr\'esentant $(P_1, \alpha,P_2)$ de $x$. Supposons d'abord
$(P_1,\alpha,P_2)\cong (P'_1, \alpha', P'_2)$ dans la cat\'egorie ${\cal
C}$. Il existe dans ce cas un couple d'isomorphismes $f_i : P_i \cong
P'_i$ tel que $\alpha'\circ(nf_1) = (nf_2) \circ \alpha$. Choisissons des
syst\`emes de coordonn\'ees ${\cal S'}_i$ sur $P'_i$. Posons
\begin{eqnarray*}
&&c = (\tr \alpha^{-1} d_{nc}\alpha, \rg(P_1,{\cal S}_1) - \rg(P_2, 
{\cal S}_2))
\text{ et } \\
&& c' =  (\tr \alpha'{}^{-1} d_{nc}\alpha', \rg(P'_1,{\cal S'}_1) - 
\rg(P'_2, {\cal
S'}_2))\,.
\end{eqnarray*}
  On a  $\partial c = \partial c' = 0$.
Introduisons
\begin{eqnarray*}
\omega_1 & = & \tr(f^{-1}_2 d_{nc}f_2 - f_1^{-1} d_{nc}f_1) \text{ et } \\
\omega_2 & = & \tr((\alpha')^{-1}d_{nc}(nf_2) d_{nc}(\alpha\circ 
nf_1^{-1}) + (nf_1)
\alpha^{-1}d_{nc}\alpha d_{nc}(nf_1^{-1})) ;
\end{eqnarray*}
on a $\omega_i \in \Omega^i(A)$ et $\partial (\omega_2,\omega_1) =
c'-c$, c'est-\`a-dire $[c'] = [c]$ dans $HH_1(A, {\bf Z}/n)$.

Le m\^eme argument montre que la classe de $c$ est ind\'ependante du choix des
syst\`emes de coordonn\'es ${\cal S}_1$ et ${\cal S}_2$ retenus sur 
$P_1$ et $P_2$.

\noindent Montrons \`a pr\'esent qu'on a un morphisme de groupes
$$D:K_0({\cal C}) \to HH_1(A; {\bf Z}/n)$$ en posant
$D((P_1,\alpha,P_2)^{\textstyle\protect\widehat{}}\ ) =
[\tr(\alpha^{-1}d_{nc}\alpha), \rg(P_1, {\cal S}_1)-\rg(P_2, {\cal 
S}_2)]$. Pour
$i=1,2,3$, soient  $t_i = (P_i,\alpha_i, Q_i)^{\textstyle\protect\widehat{}} \
$
  trois \'el\'ements de $K_0({\cal C})$ tels que
$t_1+t_2=t_3$. Apr\`es avoir fix\'e des syst\`emes de coordonn\'ees
${\cal S}_i$ et ${\cal S'}_i$ sur $P_i$ et $Q_i$ $(i=1,2)$, on choisit
${\cal S}_1\cup {\cal S}_2$ (resp. ${\cal S'}_1\cup {\cal S'}_2$) comme
syst\`eme de coordonn\'ees de $P_3$
  (resp. $Q_3$).

\noindent Posons $c_i =(\tr(\alpha^{-1}_id_{nc}\alpha_i), \rg(P_i, {\cal
S}_i)-\rg(Q_i, {\cal S}_i))$ pour $i=1,2$.

\noindent De $\alpha_3 = \alpha_1 \oplus \alpha_2$  et du choix
propos\'e pour les syst\`emes de coordonn\'ees, on tire imm\'ediatement
$c_3 = c_1 + c_2$. Cette \'egalit\'e entre cocycles entra\^\i ne$[c_1] +
[c_2] = [c_3]$.

Il reste \`a montrer $D(N) = 0$, o\`u $N$ est le sous-groupe de $K_0({\cal
C})$ tel que $K_1(A; {\bf Z}/n) = K_0({\cal C})/N$. Soient  $z_1 =
(P_1,\alpha,P_2)^{\textstyle\protect\widehat{}}\ $ et $z_2 = (P_2, \beta,
P_3)^{\textstyle\protect\widehat{}}\ $ deux
\'el\'ements de $K_0({\cal C})$;  on pose $z = (P_1, \beta\alpha,
P_3)^{\textstyle\protect\widehat{}}\ $. Apr\`es avoir choisi des syst\`emes
de coordonn\'ees
${\cal S}_i$ sur $P_i$, on pose
\begin{eqnarray*}
c_1 & = & (tr \alpha^{-1}d_{nc}\alpha, \rg(P_1, {\cal S}_1)- \rg(P_2, {\cal
S}_2)) \,, \\
c_2 & = & (\tr \beta^{-1} d_{nc}\beta, \rg(P_2, {\cal S}_2) - \rg(P_3, {\cal
S}_3)\,, \\
c & = & (\tr(\beta\alpha)^{-1} d_{nc}(\beta\alpha), \rg(P_1, {\cal
S}_1)-\rg(P_3, {\cal S}_3))\,, \\
\theta_2 & = & -\tr(\alpha^{-1} \beta^{-1}d_{nc}\beta d_{nc}\alpha) \,, \\
\theta_1 & = & 0\,.
\end{eqnarray*}
On a $c_1+c_2-c_3 = \partial (\theta_2,\theta_1)$, ce qui montre que
l'application $D$ factorise en un morphisme $ D_1^{(n)} : K_1(A, {\bf Z}/n)
\to HH_1(A;{\bf Z}/n)$. En remarquant que\break $\tr(\alpha\oplus\id)^{-1}
d(\alpha\oplus\id) = \tr (\alpha^{-1}d \alpha)$, il est imm\'ediat de 
v\'erifier
que les diagrammes commutent.
\endpreuve

\begin{cor}
Posons $\widetilde {HH_1}(A;{\bf Z}/n) = HH_1(A; {\bf
Z}/n)/\mathrm{Im}(\rho\circ D_1)$. La ``classe caract\'eristique secondaire''
$$
\overline d_1 : K_0(A)_{(n)} \to \widetilde {HH_1}(A;{\bf Z}/n)\, .
$$
d\'efinie pour $x=\partial (y)\in K_0(A)_{(n)}$
par
$\overline d_1(x)=\overline D_1(y) \ \mod \ \hbox{Im}\rho\circ D_1$ 
est un morphisme de groupes
ab\'eliens.
\end{cor}
Nous montrerons plus loin que cette classe caract\'eristique secondaire n'est
pas triviale en g\'en\'eral.

\subsection{Le cas des alg\`ebres commutatives.}

Si l'alg\`ebre $A$ est commutative, on peut
simplifier notablement la cons\-truction
de la trace de Dennis \`a coefficients en transitant par le $A$-module
$\Omega^1_{dR}(A)$ des diff\'erentielles de K\"ahler de $A$.  Fixons 
pour cela  un peu
de voca\-bulaire.

Le module des diff\'erentielles
de K\"ahler de $A$ est le $A$-module $\Omega^1_{dR}(A):=\ker \mu/(\ker \mu)^2$,
o\`u $\mu:A\otimes_kA\rightarrow A$ est la multiplication de $A$.
Ce $A$-module est
engendr\'e par
les symboles $da$, $a \in A$, soumis aux relations
$d\lambda =0$, $\lambda\in k$,
$d(a_0a_1) = a_0 da_1 + a_1 da_0$.

On sait que lorsque $A$ est commutative, la trace de Dennis
$D_1$ est essentiellement la d\'eriv\'ee logarithmique. En effet, de la
commutativit\'e de $A$,
on d\'eduit
$$HH_1(A)=\Omega_{nc}^1(A)/[A,\Omega_{nc}^1(A)].$$
L'appli\-ca\-tion $\gamma :HH_1(A)\rightarrow \Omega^1_{dR}(A)$ d\'efinie par
$\gamma([a_0d_{nc}a_1])=a_0da_1$ est
un isomorphisme.
On en d\'eduit le morphisme surjectif de $A$-bimodules
$\rho:\Omega^1_{nc}(A)\to \Omega^1_{dR}(A)$ d\'efini par
$\rho(a_0d_{nc}a_1)=a_0da_1$.

Pour $x=[P,\alpha]\in K_1(A)$, on obtient
facilement $\gamma\circ D_1(x)=\break\det(\alpha)^{-1}d(\det(\alpha)).$
En \'ecrivant $K_1(A)=A^{\times}\oplus SK_1(A)$, on en d\'eduit que 
la restriction de $D_1$
\`a $SK_1(A)$ est nulle et que la restriction de la trace de Dennis \`a
$A^\times$ est la d\'eriv\'ee logarithmique, c'est-\`a-dire que
pour $u\in
A^{\times}$, on a
$\gamma\circ D_1(u)=u^{-1}du$.

Soit $P$ un $A$-module projectif de type fini sur $A$ et soit
${\cal S}=(x_j,\varphi_j)_{1\leq j\leq r}$ un syst\`eme de 
coordonn\'ees sur $P$. La
connexion de Levi-Civita commutative de $P$ est l'application
$d'_P:P\to P\otimes_A\Omega^1_{dR}(A)$ d\'efinie pour
$x=\sum_{j=1}^rx_j\varphi_j(x)$ par
$d'_p(x)=\sum_{j=1}^rx_jd\varphi_j(x)$.
Remarquons que si $d_P:P\to P\otimes_A\Omega^1_{nc}(A)$ d\'esigne la
connexion  de Levi-Civita non commutative, le diagramme suivant est
commutatif.
$$
\xymatrix{
P\ar[r]^-{d'_P}\ar[d]^-{d_P}& P\otimes_A\Omega_{dR}^1(A)\\
P\otimes_A\Omega^1_{nc}(A)\ar[ur]_-{\id\otimes\rho}
}
$$

Soient $P$ et $Q$ deux $A$-modules projectifs de type fini et soit 
$\alpha:P\rightarrow
Q$ une  application $A$-lin\'eaire.
Soient $$\nabla:P\rightarrow P\otimes_{A}\Omega^1_{dR}(A)$$ et
$$\nabla ':Q\rightarrow Q\otimes_{A}\Omega^1_{dR}(A)$$ des connexions 
sur $P$ et $Q$
respectivement.  On d\'efinit l'application $A$-lin\'eaire
$d\alpha=d(\alpha,\nabla,\nabla ')$  de source
$P$, de but $Q\otimes_{A}\Omega^1_{dR}(A)$
par
$$d(\alpha,\nabla,\nabla ')=\nabla '\circ\alpha-(\alpha\otimes 
id)\circ \nabla.$$

Soit $\alpha :P\rightarrow Q$ un isomorphisme de $A$-module. En 
choisissant des syst\`emes de
coordonn\'ees sur $P$ et $Q$, l'application $A$-lin\'eaire
$$\alpha^{-1}d\alpha:=(\alpha^{-1}\otimes\id)\circ d(\alpha,\nabla,\nabla')$$
admet une matrice carr\'ee \`a coefficients dans $\Omega^1_{dR}(A)$ 
dont la trace
est not\'ee $tr(\alpha^{-1}d\alpha)$.

\begin{theo}
Soient $k$ un anneau commutatif unitaire,
$A$ une $k$-alg\`ebre commutative et soit  $n\geq 2$ un entier.

Soit
$$D_1^{(n)}:K_1(A;{\bf Z}/n)\rightarrow \Omega^1_{dR}(A)/(n)$$
l'application
d\'efinie pour $x=[P,\alpha ,Q]$ dans $K_1(A,{\bf Z}/n)$ par
$$D_1^{(n)}(x)=
tr(\alpha^{-1}\circ d(\alpha,n\nabla ,n\nabla '))
\ \mod\ n\Omega_{dR}^1(A),$$
o\`u $\nabla$ et $\nabla '$ sont des connexions sur $P$ et $Q$ respectivement.

Alors l'application $D_1^{(n)}$ est un morphisme de groupes ab\'eliens.

\end{theo}
\rem
Supposons de plus que l'alg\`ebre  $A$ soit int\`egre et que $n$  soit
un entier premier \`a la caract\'eristique de $A$. La suite exacte longue 1.3
conduit
\`a  l'isomorphisme $HH_1(A;{\bf Z}/n)\cong HH_1(A)/(n)$. Notons
$\overline \gamma:HH_1(A)/(n)\rightarrow \Omega^1_{dR}(A)/(n)$ l'isomorphime
induit de
$\gamma:HH_1(A)\simeq \Omega^1_{dR}(A)$.
A l'aide des connections de Levi-Civita, on v\'erifie que 
l'application $D_1^{(n)}$
s'ins\`ere dans le  diagramme commutatif
$$
\xymatrix{
  K_1(A;{\bf Z}/n)\ar[dr]_-{D^{(n)}_1}\ar[r]^-{\overline D_1} &HH_1(A;{\bf
Z}/n)\ar[d]^-{\overline\gamma}\\
  &\Omega^1_{dR}(A)/(n)
  }
$$
%

\noindent Preuve du th\'eor\`eme 23:  avec les notations de 1.1, on a
$K_1(A;{\bf Z}/n)=K_0({\cal C})/N$. Soit
$(P,\alpha, Q)$ un objet de ${\cal C}$.

Si $\nabla$ et $\nabla_1$ sont deux connexions sur $P$ et si
$\nabla '$ est une connexion sur $Q$, on a
$$\alpha^{-1}\circ d(\alpha,n\nabla,n\nabla ')-
\alpha^{-1}\circ d(\alpha,n\nabla_1,n\nabla ')
=n(\nabla-\nabla_1),$$
application $A$-lin\'eaire dont la trace est congrue \`a $0$ modulo
$n\Omega^1_{dR}(A).$

Si $\nabla$ est une connexion sur $P$ et si
$\nabla '$ et $\nabla '_1$ sont deux connexions sur $Q$, on a
$$\alpha^{-1}\circ d(\alpha,n\nabla,n\nabla ')-
\alpha^{-1}\circ d(\alpha, n\nabla,n\nabla '_1)=
\alpha^{-1}\circ n(\nabla '-\nabla '_1)\circ\alpha,$$
  application $A$-lin\'eaire dont la trace, \'egale \`a celle de
$n(\nabla '-\nabla '_1)$, est bien congrue  \`a $0$ modulo
$n\Omega^1_{dR}(A)$.

Ces deux remarques montrent que pour
  $(P,\alpha ,Q)\in Ob({\cal C})$, le choix des con\-nexions sur $P$ ou $Q$
n'intervient pas pour la d\'efinition de
$tr(\alpha^{-1}d\alpha)$ mo\-du\-lo $n\Omega^1_{dR}(A)$.
C'est pourquoi, pour all\'eger, nous omettons \`a pr\'esent de pr\'eciser
les connexions choisies.

Si les objets $(P,\alpha,Q)$ et $(P_1,\alpha_1,Q_1)$ sont
isomorphes dans la cat\'egorie ${\cal C}$, il existe des applications
$A$-lin\'eaires $f$ et $g$ telles que
$\alpha_1=ng\circ \alpha\circ nf^{-1}$.
Un rapide calcul donne
$$\tr(\alpha_1^{-1}d\alpha_1)=
n\tr(g^{-1}dg)+\tr(\alpha^{-1}d\alpha)+n\tr(fdf^{-1})$$
soit
$$\tr(\alpha_1^{-1}d\alpha_1)\equiv
\tr(\alpha^{-1}d\alpha)\  \mod\ n\Omega^1_{dR}(A),$$
ce qui montre que seule la classe d'isomorphie  de l'objet $(P,\alpha,Q)$
intervient pour la d\'efinition de la trace \`a coefficients.
Enfin,  les relations banales
$$\tr\left(
(\alpha_1\oplus\alpha_2)^{-1}d(\alpha_1\oplus\alpha_2)\right)=
\tr(\alpha_1^{-1}d\alpha_1)+\tr(\alpha_2^{-1}d\alpha_2)$$
et
$$\tr\left(
(\alpha\beta)^{-1}d(\alpha\beta)\right)=
\tr(\alpha^{-1}d\alpha)+\tr(\beta^{-1}d\beta)$$
montrent qu'on a  un morphisme de groupes
$D_1^{(n)}$ de source
$K_1(A;{\bf Z}/n)$ de but
$\Omega_{dR}^1(A)/(n)$ en posant
$$D_1^{(n)}\left([P,\alpha,Q]\right)=
\tr(\alpha^{-1}d\alpha)\ \mod \ n\Omega^1_{dR}(A).$$
\medskip

\noindent
{\bf Exemple.}
Soit $A$ un anneau commutatif. On suppose $K_0(A)={\bf Z}$.  Alors,
$K_1(A;{\bf Z}/n)=K_1(A)/(n)=A^\times/(n)\oplus SK_1(A)/(n)$.
La restriction de la trace de Dennis $D_1^{(n)}$ au facteur 
$SK_1(A)/(n)$ est nulle
tandis que la restriction de la trace de Dennis $D_1^{(n)}$ au facteur
$A^\times/(n)$ est donn\'ee par
$$D_1^{(n)}([a])=a^{-1}da\ \mod\ n\Omega^1_{dR}(A).$$
Cette situation s'applique en particulier lorsque $A$ est local.
\bigskip

Pour tout anneau commutatif, l'image de $K_1(A)$ par la trace de 
Dennis $D_1$ est
le sous-groupe
$dA^\times/A^\times$ de $\Omega^1_{dR}(A)$ engendr\'e par
$\{u^{-1}du ,\ u\in A^\times\}$. Du th\'eor\`eme 23, on d\'eduit:

\begin{cor}
Soient $k$ un anneau commutatif unitaire, $A$ une $k$-alg\`ebre
commutative et $n$ un entier. On  d\'esigne par
$A^\times$ le groupe des  unit\'es de $A$ et par
$dA^\times/A^\times$ le sous-groupe de
$\Omega^1_{dR}(A)$ engendr\'e par
$\{u^{-1}du, \ u\in A^\times\}$. Soit
$S$ le sous-groupe de $\Omega^1_{dR}(A)$ engendr\'e par $n\Omega^1_{dR}(A)$ et
$dA^\times/A^\times$. La ``classe caract\'eristique secondaire''
$$  d_1^{(n)}:
\tilde K_0(A)_{(n)}\rightarrow
\Omega^1_{dR}(A)/S$$
d\'efinie pour $x=\partial(y)\in\tilde K_0(A)_{(n)}$
par
$d_1^{(n)}(x)=D_1^{(n)}(y)\ \mod\ S$
est un morphisme de groupes ab\'eliens.
\end{cor}
Nous verrons plus bas que ce morphisme n'est pas trivial.
\medskip

\subsection{Les traces d'ordre sup\'erieur.}

Soit $A$ une $k$-alg\`ebre. D\'esignons par $HH_*(A)$, $HC_*^-(A)$ et
$HC_*^{per}(A)$ res\-pec\-ti\-vement  les homologies de Hochschild,
``cyclique n\'egative'' et ``cyclique p\'eriodique'' de $A$ (voir 
[17], 1.1 et 5.1
pour les d\'efinitions). Soit
$D_*:K_*(A)\to HH_*(A)$ la trace de Dennis [10] et
$\gamma_*:K_*(A)\to HC_*^-(A)$ le ca\-rac\-t\`ere de Chern universel
([8], II). Pour tout
$r\geq 0$, Weibel ([25], p. 541) a montr\'e qu'on a un diagramme commutatif
$$
\xymatrix{
K_r(A)\ar[r]^-{\gamma_r}\ar[d]^-{D_r} &HC_r^-(A)\ar[dl]^-{h_r}\ar[r]^-{I_r}
&HC^{per}_r(A)\\ HH_r(A) & &
  }\leqno{(0)_r}
  $$

Nous montrons ici que pour tout $r\geq 1$ et tout $n\geq 2$, il 
existe un diagramme commutatif
$$
\xymatrix{
K_r(A;{\bf Z}/n)\ar[r]^-{\gamma_r^{(n)}}\ar[d]^-{D_r^{(n)}}
&HC_r^-(A;{\bf Z}/n)\ar[dl]^-{h_r^{(n)}}\ar[r]^-{I_r^{(n)}} &HC^{per}_r(A;{\bf
Z}/n)\\ HH_r(A;{\bf Z}/n) & &
  }\leqno{(0)^{n}_r}
  $$
reliant la $K$-th\'eorie \`a coefficients de $A$ aux diverses 
homologies \`a coefficients  au moyen
d'applications $\gamma_r^{(n)}$ et $D_r^{(n)}$ qui sont pr\'ecis\'ees plus bas.

Notons $C_*(A)$, $CC_*(A)^-$ et $CC_*^{per}(A)$ les complexes de 
cha\^{\i}nes des homologies de
Hochschild, cyclique n\'egative et p\'eriodique de $A$. La 
correspon\-dance de Dold-Kan
([26], 8.4) $DK:Ch_+\to AbS$ r\'ealise une \'equivalence entre
la cat\'egorie
$Ch_+$ des complexes de cha\^{\i}nes gradu\'es positivement et la 
cat\'egorie des groupes ab\'eliens
simpliciaux. Notons $\mid\ \mid:AbS\to Top$ le foncteur r\'ealisation
g\'eom\'etrique et
$X:Ch_*\to Top$ la compos\'ee de ces deux foncteurs.
Les espaces ${\cal H}(A)=X(C_*(A))$,
${\cal H}^-(A)=X(CC_{*}^{-}(A))$ et ${\cal 
H}^{per}(A)=X(CC_{*}^{per}(A))$ sont des espaces
classifiants pour les homologies de Hochschild ou cycliques. Pour 
$r\geq 0$, on a
$\pi_r({\cal H}(A))=HH_r(A)$,
$\pi_r({\cal H}^-(A))=HC_*^-(A)$ et
$\pi_r({\cal H}^{per}(A))=HC_*^{per}(A)$. D'apr\`es [25], prop. 4.3, on a

\begin{prop}
Pour toute alg\`ebre $A$, il existe dans la cat\'egorie $Top$ un 
diagramme commutatif
$$
\xymatrix{
BGLA^+\ar[r]^-{\gamma}\ar[d]^-{D} &{\cal H}^-(A)\ar[dl]^-{h} \ar[r]^-{I}
&{\cal H}^{per}(A)\\ {\cal H}(A) & &
  }\leqno{(1)}
  $$
tel que pour $r\geq 0$, en appliquant le  foncteur $\pi_r(-)$ au 
diagramme $(1)$,
on obtienne le diagramme $(0)_r$.
\end{prop}

Il est naturel de chercher \`a obtenir le diagramme $(0)^n_r$ par 
application du foncteur
$\pi_r(-;{\bf Z}/n)$ au diagramme $(1)$. Ce foncteur n'est d\'efini 
que pour $r\geq 2$.
Une d\'esuspension permet  l'\'etude du cas $r=1$.

Rappelons (voir [19]), que pour $r\geq 2$, on pose
$\pi_r(-;{\bf Z}/n)=[M^r_n,-]$
o\`u $M^r_n$ est l'espace de Moore $S^{r-1}\cup_\alpha e^r$ avec 
$\alpha$ de degr\'e $n$.

Pour $r\geq 2$, on pose $K_r(A;{\bf Z}/n)=\pi_r(BGLA^+;{\bf Z}/n)$. 
La  structure de $H$-espace de
$BGLA^+$ permet de d\'efinir une application $.n:BGLA^+\to BGLA^+$ 
induisant la multiplication par
$n$ sur
$\pi_r(BGLA^+)$. La fibre homotopique ${\cal F}$ de $.n:BGLA^+\to 
BGLA^+$ est telle que pour $r\geq
1$, on a
$\pi_r({\cal F})=\pi_{r+1}(BGLA^+;{\bf Z}/n)$.

Rappelons \'egalement que l'homologie \`a coefficients $H_*(C_*;{\bf 
Z}/n)$ d'un
complexe de
cha\^{\i}nes $C_*$ (d\'efinie en 1.3 comme l'homologie du c\^one 
$co(.n)$ de la multiplication
$.n:C_*\to C_*$) est telle que pour $r"\geq 2$, on ait
$H_r(C_*;{\bf Z}/n)=\pi_r(X(C_*);{\bf Z}/n)$.

En conclusion, on a

\begin{prop}
Pour $r\geq 2$, le diagramme $(0)^{n}_r$ est obtenu par application du foncteur
$\pi_r(-;{\bf Z}/n)$ au diagramme $(1)$.
\end{prop}
  Pour $r=1$, contentons-nous de traiter le cas de l'homologie 
cyclique n\'egative en d\'etaillant
la construction de l'application $\gamma_1^{(n)}$ du diagramme 
$(0)^{n}_1$. Les constructions
pour l'homologie de Hochschild ou p\'eriodique sont analogues.

Soit $SA$ le c\^one de l'anneau $A$ au sens de [13], p. 269. D'apr\`es [23],
prop. 3.2, pour
$r\geq1$, on a
$K_r(SA)\cong K_{r-1}(A)$. En particulier,\break 
$\pi_2(BGL(SA)^+)\cong K_1(A)$.
D'apr\`es [27], thm. 10.1 et sa remarque 2, pour $r\geq 1$, on a
$HC_r^-(SA)\cong HC^{-}_{r-1}(A)$. En particulier, $\pi_2({\cal 
H}^-(SA))\cong HC_1^-(A)$.

L'application $\gamma':BGL(SA)^+\to {\cal H}^-(SA)$ obtenue en appliquant la
proposition 25 \`a l'anneau $SA$ permet de d\'efinir
$\gamma_1^{(n)}=\pi_1(\gamma';{\bf Z}/n)$,\break
$K_1(A;{\bf Z}/n)=\pi_2(BGL(SA)^+;{\bf Z}/n)$ et
$HC^{-}_1(A;{\bf Z}/n))=\pi_2({\cal H}^-(SA);{\bf Z}/n)$.
On a
$\gamma_{1}^{(n)}:K_1(A;{\bf Z}/n)\to HC_1^{-}(A;{\bf Z}/n).$
Les suites exactes longues de Barratt pour l'homotopie \`a coefficients
([19], p. 3) nous donnent

\begin{prop}
On a le diagramme commutatif naturel \`a lignes exactes

\newsavebox{\maboitecinq}
\savebox{\maboitecinq}(350,80)[]{\scriptsize
$
\xymatrix{
K_1(A)\ar[d]^-{\gamma_1} \ar[r]^-{.n} &
K_1(A) \ar[d]^-{\gamma_1} \ar[r]^-\rho &
K_1(A;{\bf Z}/n) \ar[d]^-{\gamma^{(n)}_1} \ar[r]^-{\partial} &
  K_0(A)\ar[d]^-{\gamma_0}\ar[r]^-{.n} &
K_0(A)\ar[d]^-{\gamma_0} \ar[r] &
K_0(A)/(n)\ar[d]^-{\overline \gamma_0}\\
HC_1^-(A)\ar[r]^-{.n} &
HC_1^-(A)\ar[r]^-{\rho} & HC_1^{-}(A;{\bf Z}/n) \ar[r]^-{\partial}
&HC_0^-(A) \ar[r]^-{.n} &HC_0^-(A)\ar[r] & HC_0^-(A;{\bf Z}/n)
  }
$}
\usebox{\maboitecinq}
\end{prop}
\begin{rem}
Il est raisonnable de conjecturer que l'application compos\'ee
  $$\xymatrix{
K_1(A;{\bf Z}/n)\ar[r]^-{\gamma_1^{(n)}}&
  HC_1^-(A;{\bf Z}/n)\ar[r]^-{h_1^{(n)}}&
  HH_1(A;{\bf Z}/n) }$$
est celle d\'ecrite en 1.5.
\end{rem}

\section{\'Etude de l'anneau des entiers d'un corps de nombres.}

Soit $A$ l'anneau des entiers d'un corps de nombres $F$.
D'apr\`es [3], on a $K_1(A)=A^\times$. L'extension $(\ddagger)$ 1.1 s'\'ecrit
donc
  $$
\xymatrix{
1 \ar[r]  &A^\times/ (n) \ar[r]^-{\rho} & K_1(A;{\bf Z}/n)
\ar[r]^-{\partial}  & Cl(A)_{(n)}\ar[r] & 0\,
}\leqno{(\dagger )}
$$
ce qui montre que $K_1(A;{\bf Z}/n)$ est une extension de la 
$n$-torsion du groupe des classes de
$A$ par un quotient du groupe des unit\'es de $A$.

La trace de Dennis $D_1^{(n)}$ et la classe
carac\-t\'e\-ristique secondaire
$d_1^{(n)}$ cons\-truites en 1.6 s'ins\`erent dans le diagramme commutatif
$$
\xymatrix{
  K_1(A;{\bf Z}/n) \ar[d]_-{D_1^{(n)}}\ar[r]^-{\partial} & Cl(A)_{(n)} 
\ar[d]^-{ d_1^{(n)}}
\ar[r]  & 0  \\
\Omega^1_{dR}(A)/(n)  \ar[r] & \Omega^1_{dR}(A)/S \ar[r]  &  0
}
$$

En 2.1, on montre que $K_1(A;{\bf Z}/n)$ est isomorphe \`a un groupe 
not\'e\break
${\cal U}(A;{\bf Z}/n)$ construit \`a partir d'id\'eaux fractionnaires de $A$.
  Le lemme ($N$-$N_1$) de la section 2.2 fournit un crit\`ere de construction
d'\'el\'ements de
${\cal U}(A;{\bf Z}/n)$. En 2.3, on d\'ecrit le lien entre
$K_1(A;{\bf Z}/n)$ et le groupe des ad\`eles restreints $\widehat A$ 
de $A$. Les
traces $D_1^{(n)}$ et $d_1^{(n)}$ sont d\'etaill\'ees en 2.5.

\subsection{Description de $K_1(A;{\bf Z}/n)$ en termes d'id\'eaux.}

Soit $A$ un anneau de Dedekind, de corps des fractions $F$ et soit
$n\geq 2$ un entier. On d\'esigne par $I(A)$ le mono\"{\i}de des id\'eaux
fractionnaires de $A$.  Pour
$x\in F$, on pose
$[x]=x\
\mod
\ F^{\times(n)}$. Consid\'erons le sous-groupe ${\cal U}(A;{\bf Z}/n)$ de
$F^\times/(n)$ d\'efini par
$${\cal U}(A;{\bf Z}/n):=\{x\in F^\times/(n)\mid\exists\ I\in I(A),\
xA=I^n\}\cdot$$

\begin{theo}
Soit $A$ un anneau de Dedekind.
Alors, on a un isomorphisme
$$K_1(A;{\bf Z}/n)\cong{\cal U}(A;{\bf Z}/n)\oplus SK_1(A)/(n).$$
\end{theo}

\preuve
on sait que dans un anneau de Dedekind, tout module $P$ projectif de 
type fini et de
rang $r$ est de la forme $P=(r-1)A\oplus I$ o\`u $I$ est un id\'eal
fractionnaire de
$A$. Soit $[P,\alpha, L]$ un \'el\'ement de $K_1(A;{\bf Z}/n)$ avec
$P=(r-1)A\oplus I$ et $L=rA$. L'isomorphisme $\alpha:nP\cong nL$ montre
qu'il exis\-te $x\in F^\times$ tel que $I^n=xA$. D\'efinissons
$\det^{(n)}:K_1(A;{\bf Z}/n)\to{\cal U}(A;{\bf Z}/n)$ par
$\det^{(n)}([P,\alpha,L])=[x].$
On v\'erifie que $\det^{(n)}$ est bien d\'efinie et que c'est un morphisme de
groupes. Pour obtenir le facteur direct
${\cal U}(A;{\bf Z}/n)$ de $K_1(A;{\bf Z}/n)$, on construit un morphisme
$s:{\cal U}(A;{\bf Z}/n)\to K_1(A;{\bf Z}/n)$ tel que
$\det^{(n)}\circ s=\id_{{\cal U}(A;{\bf Z}/n)}$. Soit
$[x]\in{\cal U}(A;{\bf Z}/n)$ avec $xA=I^n$ o\`u $I$ est un id\'eal 
fractionnaire de
$A$. On pose $s([x])=[A,\id_{(n-1)A}\oplus \cdot x,I]$. On v\'erifie 
que $s$ est bien
d\'efinie. Pour montrer $s([x][y])=s([x])+s([y])$, posons
$xA=I^n$, $yA=J^n$, $s([x])=u$, $s([y])=v$,
$s([x][y])=w$, avec
$u=[A,f,I]$, $f=\id_{(n-1)A}\oplus \cdot x$,
$v=[A,g,J]$, $g=\cdot y\oplus\id_{(n-1)A}$,
$w=[A,h,IJ]$, $h=\cdot xy\oplus\id_{(n-1)A}$.
On a $u+v=[A\oplus A, h_1, I\oplus J]$ o\`u $h_1$ est l'isomorphisme
$\id_{(n-1)A}\oplus \cdot x\oplus\cdot y\oplus \id_{(n-1)A}$.
Posons $h_2=\id_{nA}\oplus h$. Dans la cat\'egorie ${\cal C}$, ({\sl 
cf.} 1.1), les
objets
$(A\oplus A,h_1,I\oplus J)$ et
$(A\oplus A,\id_{nA}\oplus h, A\oplus IJ)$
sont isomorphes.
Dans $K_1(A;{\bf Z}/n)$, on a par cons\'equent
$u+v=[A\oplus A,h_1,I\oplus J]=[A,\id_{nA},A]+[A,h,IJ]=0+w=w$, ce qui montre
que $s$ est un morphisme de groupes. On v\'erifie
la relation $\det^{(n)}\circ s=\id_{{\cal U}(A;{\bf Z}/n)}$, ce qui montre

\smallskip
\centerline{$K_1(A;{\bf Z}/n)={\cal U}(A;{\bf Z}/n)\oplus\ker\ \det^{(n)}.$}
Introduisons le diagramme
$$\xymatrix{
K_1(A)\ar[r]^-{\cdot n}\ar[d]^-{\det}&
K_1(A)\ar[r]^-{\rho}\ar[d]^-{\det} &
K_1(A;{\bf Z}/n)\ar[r]^-{\partial}\ar[d]^-{\det^{(n)}}&
K_0(A)\ar[r]^-{\cdot n}\ar[d]^-{\id}&
K_0(A)\ar[d]^-{\id}\\
A^\times\ar[r]^-{\cdot n}&
A^\times \ar[r]^-{\rho'}&
{\cal U}(A;{\bf Z}/n)\ar[r]^-{\partial'}&
K_0(A)\ar[r]^-{\cdot n}&
K_0(A)
}$$
o\`u les applications $\rho'$ et $\partial'$ sont d\'efinies par
$\rho'(a)=[a]$, $\partial'([x])=[I^{-1}].$
Ce diagramme est commutatif \`a lignes exactes. La suite exacte des noyaux des
fl\`eches verticales conduit \`a
$\ker\hbox{det}^{(n)}\cong SK_1(A)/(n).$

\begin{cor}
Soit $A$ l'anneau des entiers d'un corps de nombres $F$. Alors
$$K_1(A;{\bf Z}/n)\cong {\cal U}(A;{\bf Z}/n).$$
  \end{cor}
En effet, d'apr\`es [BMS], on a
$SK_1(A)=0$.

\subsection{Le lemme ($N$-$N_1$) de construction d'\'el\'ements de 
$K_1(A;{\bf Z}/n)$.}

Soit $L/F$ une extension de corps de nombres de degr\'e $\ell$. On pose $A =
{\cal O}_F$ et on note
$B$ la fermeture int\'egale de $A$ dans $L$.

Pour $z \in L$, on d\'esigne par $\mu_z : L\to L$ la multiplication par $z$.
Les quantit\'es $N_j(z) \in  L$ sont d\'efinies par $\det(X  \id_L -\mu_z)
= \sum^\ell_{j=0} (-1)^{\ell-j} N_j(z) X^j$. En particulier,
$$
N_\ell(z) = 1 \,, N_{\ell-1}(z) = \tr_{L/F}(z)\,, N_0(z) = N_{L/F}(z) 
= N(z) \,.
$$
\begin{prop}
Soit $z \in  L$ (resp. $B$) et $h \in  F$ (resp. $A$). Alors on a
$$
N(z+h) = N(z) + N_1(z)h + h^2 \varepsilon (z,h) \text{ avec } \varepsilon (z,h)
\in  F (\text{ resp. } A)\,.
$$
\end{prop}

Remarquons qu'on a la formule commode
$$
N_1(z) = \left(\left(\frac{d}{dh}\right)_{h\in F} N(z+h)\right)_{h=0}\,.
$$

\begin{lemme}
Soit $L/F$ une extension de corps de nombres et soient $A$ l'anneau des entiers
de $F$ et
$B$ la fermeture int\'egrale de $A$ dans $L$. Soit $n$ un entier $\geq 2$.
Consid\'erons un
\'el\'ement
$u$ de $B$ tel que
\begin{eqnarray*}
&&N(u) = e a^n \text{ avec } e \in  A^\times \,, \, a \in  A \\
\text{ et } && (N(u), N_1(u)) = A.
\end{eqnarray*}
En d\'esignant par $[u]$ la classe de $u$ dans $L^\times/ (n)$ ,
on a alors
$$[u] \in  K_1(B;{\bf Z}/n).$$
\end{lemme}
\preuve
  L'hypoth\`ese $(N(u),N_1(u))=A$ signifie que $u$ est premier \`a
tous ses conjugu\'es. En effet, si $\sigma:L\to {\bf C}$ d\'esigne un
$F$-plongement de $L$ (c'est-\`a-dire $\sigma\mid_F=\id$), on a
$N(u)=\prod_{\sigma}\sigma(u)$,
$N_1(u)=\left(\prod_{\sigma}\sigma(u)\right)
\left(\sum_{\sigma}\sigma(u)^{-1}\right)$ et
$P_u(X)=\prod_\sigma(X-\sigma(u))$.
Soit $\mathfrak{p}$ un id\'eal premier de $A$. On a
$\mathfrak{p}\mid (N(u),N_1(u))$
si et seulement
si $P_u(X)\cong X^2Q(X)\ \mod\   \mathfrak{p}$.
Ceci signifie qu'il existe deux plongements $\sigma_1$ et $\sigma_2$ tels
que
$(\sigma_1(u),\sigma_2(u))\subset \mathfrak{p}$. En posant
$\tau=\sigma_1^{-1}\sigma_2$ et
$\mathfrak{q}=\sigma_1^{-1}(\mathfrak{p})$, on en d\'eduit
$(u,\tau(u))\subset\mathfrak{q}$.

Pour montrer que $[u]$ appartient \`a $K_1(A;{\bf Z}/n)$, on remarque que
la puissance $n$-i\`eme de  l'id\'eal fractionnaire $I=(u,a)$ de
$B$ est principale. En effet
$I^n=(u^n,N(u))=(u^n,u\prod_{\sigma\not=\id}\sigma(u))$; et puisque
$(u,\sigma(u))=B$, on en d\'eduit $I^n=uB$. Ceci montre
$[u]\in {\cal U}(B;{\bf Z}/n)$.
\endpreuve

\subsection{Description de  $K_1(A;{\bf Z}/n)$ en termes d'ad\`eles.}

Soit $A$ l'anneau des entiers d'un corps de nombres $F$. Notons $Spec(A)$ le
spectre premier de $A$. Pour
$\mathfrak{p}$ dans
$Spec(A)$, on d\'esigne par $\widehat A_{\mathfrak{p}}$ le compl\'et\'e
$\mathfrak{p}$-adique de l'anneau de valuation discr\`ete $A_\mathfrak{p}$.
L'anneau $\widehat A=\prod_{\mathfrak{p}\in Spec(A)}\widehat A_\mathfrak{p}$,
appel\'e anneau des ad\`eles restreints de $A$, s'ins\`ere dans le diagramme
commutatif
$$\xymatrix{
A\ar@{^{(}->}[r]^-{}\ar@{^{(}->}[d]_-{}&
{\widehat A}\ar@{^{(}->}[d]^-{}\\
F\ar[r]_-{}& {\widehat F}
}$$
o\`u on a pos\'e $\widehat F=F\otimes_A\widehat A$.

Consid\'erons le diagramme
$$\xymatrix{
{\cal U}(A;{\bf Z}/n)\ar[r]^-{j}\ar[d]_-{\imath}&
F^\times/(n)\ar[d]^-{\overline\imath}\\
\widehat A^\times/(n)\ar[r]_-{\overline j}&
\widehat F^\times/(n)
}$$
Les applications
$\imath$, $j$ et $\overline j$ sont trivialement injectives. Dans l'anneau
$\widehat A$, on s'est restreint aux places archim\'ediennes. D'apr\`es [1],
Chap. X.I, dans cette situation, l'application $\overline \imath$ est 
\'egalement
injective. La somme amalgam\'ee
$\displaystyle\widehat A^\times/(n)\oplus_{\widehat
F^\times/(n)}F^\times/(n)$ est donc \'egale \`a
$\widehat A^\times/(n)\cap F^\times/(n)$. Soit $x\in F^\times/(n)$. 
L'\'el\'ement
$[x]=x\ \mod\ F^{\times(n)}$ de $F^\times/(n)$ appartient \`a
$\widehat A^\times/(n)$ si et seulement si $n\mid
v_\mathfrak{p}(x_\mathfrak{p})$ pour tout $\mathfrak{p}\in Spec(A)$,
c'est-\`a-dire qu'on a
$xA=I^n$ avec $I$ id\'eal fractionnaire. En conclusion, nous avons:
\begin{prop}
Soit $A$ l'anneau des entiers d'un corps de nombres $F$ et soit $n\geq 2$ un
entier. Alors le groupe $K_1(A;{\bf Z}/n)$ s'identifie au sous-groupe
$\widehat A^\times/(n)\cap F^\times/(n)$ de $\widehat F^\times/(n)$.
\end{prop}

\medskip
\begin{cor}
Soit $A$ l'anneau des entiers d'un corps de nombres $F$. Posons $\widehat
A=\prod_\mathfrak{p}\widehat A_\mathfrak{p}$. Alors l'application
$K_1(A;{\bf Z}n)\to K_1(\widehat A;{\bf Z}/n)$ induite par
$A\to \widehat A$ est l'inclusion
$\widehat A^\times/(n)\cap F^\times/(n)\to \widehat A^\times/(n)$.
\end{cor}
\preuve Calculons
$K_1(\widehat A;{\bf Z}/n)$.
Rappelons pour cela que si $(A_i)_{i\in I}$ est une famille d'anneaux 
commutatifs
de rang stable $d\geq 2$ au sens de [2], p. 231, on a
$K_1(\prod_{i\in I}A_i)\cong\prod_{i\in I}K_1(A_i)$ et
$\widetilde K_0(\prod_{i\in I}A_i)\cong \prod_{i\in I}\widetilde K_0(A_i).$
Les anneaux de la famille
$(\widehat A_\mathfrak{p})_{\mathfrak{p}\in Spec(A)}$ sont tous de rang stable
$d=2$. De
$\widetilde K_0(\widehat A_\mathfrak{p})_{(n)}=0$ et
$K_1(\widehat A_\mathfrak{p})=\widehat A_\mathfrak{p}^\times$, on tire
$K_1(\widehat A)=\widehat A^\times$.
L'extension ($\ddagger$ 1.1) nous m\`ene \`a
$$K_1(\widehat A;{\bf Z}/p)=\widehat A^\times/(n).$$

\subsection{Description de $\Omega^1_{dR}(A)/(n).$}

On d\'esigne toujours par $A$ l'anneau des entiers d'un corps de nombres $F$.
L'anneau $A\otimes {\bf Z}_p$ est toujours consid\'er\'e comme une
alg\`ebre sur l'anneau ${\bf Z}_p$ des entiers $p$-adiques.
D\'esignons par $\delta$ le discriminant du corps $F$.
D'apr\`es [15], prop. 1.5, on a les \'egalit\'es suivantes
$$\Omega^1_{dR}(A)=\Omega^1_{dR}(\widehat A)=\oplus_{p\mid
\delta}\Omega^1_{dR}(A\otimes{\bf Z}_p).$$

Supposons $n$ et $p$  premiers entre eux,  alors $n$ appartient \`a $({\bf
Z}_p)^\times$ et
$\Omega^1_{dR}(A\otimes{\bf Z}_p)/(n)=0$. On en d\'eduit:

\begin{prop}
Soit $F$ un corps de nombres  d'anneaux d'entiers $A$ et de 
discriminant $\delta$.
Soit $n>1$ un entier.

\noindent a) Si $n\mid \delta$, alors
$\Omega^1_{dR}(\widehat A)/(n)=\Omega^1_{dR}(A)/(n)=\oplus_{p\mid
(n,\delta)}\Omega^1_{dR}(A\otimes{\bf Z}_p)/(n)$.

\noindent En particulier, si $p$ est un nombre premier ramifi\'e dans $A$, on a
$$\Omega^1_{dR}(A)/(p)\cong
\Omega^1_{dR}(A\otimes{\bf Z}_p)/(p).$$

\noindent b) si $(n,\delta)=1$, alors
$\Omega^1_{dR}(\widehat A)/(n)=\Omega^1_{dR}(A)/(n)=0.$

\noindent En particulier, si $p$ est un nombre premier non ramifi\'e 
dans $A$, alors on a
$\Omega^1_{dR}(A)/(p)=0.$
\end{prop}

\subsection{Description de la trace de Dennis \`a coefficients}

Soit $A$ l'anneau des entiers d'un corps de nombres $F$ et soit $n$ un diviseur
du discriminant de $F$. Pour les \'el\'ements de $K_1(A;{\bf Z}/n)$ obtenus
gr\^ace au lemme ($N$-$N_1$), il est facile de d\'ecrire la trace
de  Dennis \`a
coefficients. Supposons que $u\in A$ satisfasse aux hypoth\`eses du lemme
($N$-$N_1$) et que de plus $N(u)$ soit un entier premier \`a $n$. Posons
$v=\prod_{\sigma\not=\id}\sigma(u)$. On a $v\in A$ et $uv=N(u)$. Dans
$\Omega^1_{dR}(A)/(n)$, on en d\'eduit
$D_1^{(p)}([u])=N(u)^{-1}vdu\ \mod\ n$.

\medskip

L'\'egalit\'e $K_1(A;{\bf Z}/n)\cong F^\times/(n)\cap\widehat A^\times/(n)$
permet \'egalement de d\'ecrire localement la trace de Dennis  \`a 
coefficients.
Pour cela, on remarque que
le diagramme suivant est commutatif.
$$
\xymatrix{
K_1(A;{\bf Z}/n) \ar@{^{(}->}[d]\ar[r]^{D_1^{(n)}} &
\Omega^1_{dR}(A)/(n)\ar@{=}[d]\\
K_1({\widehat A};{\bf Z}/n) \ar[r]^{D_1^{(n)}}
&\Omega^1_{dR}({\widehat A})/(n) }
$$
  Pour conna\^{\i}tre la trace de Dennis \`a coefficients de $A$, il 
suffit donc de
conna\^{\i}tre celle de ${\widehat A}$.
Soit $p$ un nombre premier ramifi\'e dans $A$. Localement, la trace de Dennis
\`a coefficients est essentiellement une d\'eriv\'ee logarithmique modulo
$p$. En effet, on \'ecrit
$$K_1(\widehat A;{\bf Z}/p)=
\prod_{\mathfrak{p}\cap{\bf Z}\not=(p)}\widehat A^\times_\mathfrak{p}/(p)
\oplus
\prod_{\mathfrak{p}\cap{\bf Z}=(p)}\widehat A^\times_\mathfrak{p}/(p). $$
La restriction de $D_1^{(p)}$   au premier facteur de cette d\'ecomposition
est
\'evidemment nulle puisque
pour $\mathfrak{p}\cap{\bf Z}\not=(p)$, on a
$\Omega^1_{dR}(\widehat A_\mathfrak{p})/(p)=0$.
Pour $\mathfrak{p}\cap{\bf Z}=(p)$ et
$[u_\mathfrak{p}]\in \widehat A^\times_\mathfrak{p}/(p)$, d'apr\`es l'exemple
du th\'eor\`eme 23, on a
$$D^{(p)}_1([u_\mathfrak{p}])=
u_\mathfrak{p}^{-1}du_\mathfrak{p}
\ \mod\  p\Omega^1_{dR}(\widehat A_\mathfrak{p}).$$

\section{Applications aux corps de petit degr\'e.}

\subsection{Un th\'eor\`eme de Y. Yamamoto.}
L'\'egalit\'e $K_1(A;{\bf Z}/n)={\cal U}(A;{\bf Z}/n)$ et le
lemme ($N$,$N_1$)
  permettent de retrouver un
th\'eor\`eme montr\'e par  Y. Yamamoto [28] \`a l'aide de m\'ethodes 
distinctes.
\begin{theo}
Soit $F$ un corps de nombres quadratique d'anneau d'entiers $A$, de 
discriminant $\delta$
et soit $n$ un entier impair.
On suppose qu'il existe deux couples
$(\alpha,b)$ et $(\alpha',b')$ dans ${\bf Z}^2$ satisfaisant aux
relations
$$\alpha^2-4b^n={\alpha'}^2-4{b'}^n=\delta
$$
avec $(\alpha,b)=(\alpha',b')=1$.
On suppose de plus que pour tout diviseur premier $p$ de $n$, les conditions
ci-dessous sont satisfaites.
\begin{itemize}
  \item[a)] $\alpha$ (resp. $\alpha'$) n'est pas une puissance $p$-i\`eme
modulo $b$ (resp. $b'$);
\item[b)] $(\alpha+\alpha')/2$ est une puissance
$p$-i\`eme modulo $b$ et modulo $b'$.
\end{itemize}
Alors:

Si $ \delta<-4$
le groupe des classes de $A$ contient un sous groupe isomorphe \`a
${\bf Z}/n\oplus {\bf Z}/n$.

Si $\delta>0$, le groupe des classes de $A$ contient un sous groupe
isomorphe \`a ${\bf Z}/n$.
\end{theo}

\preuve L'application
$f:A\to {\bf Z}/b$ d\'efinie par
$f((x+y\sqrt{\delta})/2)=\hfill\break(x+y\alpha)/2$ est un morphisme d'anneaux.
On note
$f_1^{(d))}:K_1(A;{\bf Z}/d)\to K_1({\bf Z}/b;{\bf Z}/d)$ 
l'application induite par
$f$ en $K$-th\'eorie \`a coefficients $d$. Remarquons que
$K_1(\noindent {\bf Z}/b;{\bf Z}/d)=({\bf Z}/b)^\times/(d)$.
L'\'el\'ement $u=(\alpha+\sqrt{\delta})/2$ de $A$ est de norme $N(u)=b^n$, de
trace $\hbox{tr}(u)=N_1(u)=\alpha$. D'apr\`es le lemme ($N$-$N_1$), pour tout
diviseur $d$ de $n$, l'\'el\'ement $[u]\in F^\times/(d)$ appartient \`a
$K_1(A;{\bf Z}/d)$. Soit $p$ un diviseur premier de $n$. De 
$f(u)=\alpha$, on d\'eduit
$f_1^{(p)}([u])=[\alpha]$, quantit\'e distincte de $1$ d'apr\`es 
l'hypoth\`ese a). Pour tout
diviseur premier $p$ de $n$, l'\'el\'ement $[u]$ de $K_1(A;{\bf 
Z}/p)$ n'est donc pas
trivial. Montrons que $[u]\in F^\times /(n)$ d\'efinit un \'el\'ement 
d'ordre $n$ de
$K_1(A;{\bf Z}/n)$. Supposons $[u]$ d'ordre $m$ avec $1\leq m<n$. Il 
existe un nombre
premier $p$ tel que $mp\mid n$. De $[u]^{n/p}=1$, on tire
$u^{n/p}=z^n$ avec $z\in A^\times$, soit encore $u\in F^{\times (p)}$ 
et donc $[u]$
trivial dans
$K_1(A;{\bf Z}/p)$. On vient de montrer que ceci est impossible. On a 
donc $[u]$ d'ordre
$n$ dans $K_1(A;{\bf Z}/n)$. Le sous-groupe $H$ de $K_1(A;{\bf Z}/n)$ 
engendr\'e par
$[u]$ est donc isomorphe \`a ${\bf Z}/n$.

On introduit de mani\`ere analogue $u'=(\alpha'+\sqrt{\delta})/2$ et on obtient
de
m\^eme  un sous-groupe $H'$ de $K_1(A;{\bf Z}/n)$, \'egalement isomorphe \`a
${\bf Z}/n$. Pour montrer la somme directe $H\oplus H'$dans 
$K_1(A;{\bf Z}/n)$, on
remarque que
$f(u')=(\alpha'+\alpha)/2\in {\bf Z}/b$. Si $u'\in H$, c'est-\`a-dire 
$u'=u^m$ avec
$1\leq m<n$, l'\'egalit\'e $f_1^{(p)}([u'])=f_1^{(p)}([u])^m$, 
satisfaite pour tout
diviseur premier
$p$ de $n$, s'\'ecrit encore
$[(\alpha+\alpha')/2]=[\alpha]^m$, ce qui donne $[\alpha]^m=1$ 
d'apr\`es l'hypoth\`ese b).
On en d\'eduit comme ci-dessus qu'il existe un nombre premier $p$ tel 
que $mp\mid n$
pour lequel
$\alpha\in ({\bf Z}/b)^{\times(p)}$, ce qui fournit la contradiction 
recherch\'ee. On
montre de m\^eme $u\notin H'$. En conclusion, sous les hypoth\`eses 
propos\'ees, le
groupe
$K_1(A;{\bf Z}/n)$ contient un sous-groupe isomorphe \`a ${\bf Z}/n\oplus {\bf
Z}/n$. De l'extension
$(\dagger)$ p. 20, on d\'eduit que
si $F$ est imaginaire, $Cl(A)_{(n)}$ contient  ${\bf Z}/n\oplus {\bf 
Z}/n$ en facteur
direct.
Si $\delta>0$, $A^\times/(n)$
est isomorphe \`a ${\bf Z}/n$ et donc $Cl(A)_{n}$ contient  ${\bf 
Z}/n$ en facteur
direct.

\rem
Dans le cas o\`u $F$ est r\'eel d'unit\'e fondamentale $\varepsilon$ 
telle qu'il existe
un diviseur premier $p$ pour lequel
$f(\varepsilon)\in ({\bf Z}/b)^{\times(p)}$, le groupe des classes contient un
sous-groupe isomorphe \`a ${\bf Z}/n\oplus {\bf Z}/n$. En effet, soit 
$t=\partial([u])$
toujours avec $u=(\alpha+\sqrt{\delta})/2$. Montrons que $t$ est 
d'ordre $n$ dans
$Cl(A)$.  Supposons $t^m=0$ avec $1\leq m<n$. on en d\'eduit
$[u]^m\in \hbox{ker\ }\partial$, soit $u^m=\varepsilon^l$, ce qui donne
$[\alpha]^m=f_1^{(p)}([u]^m)=f_1^{(p)}(\varepsilon)^l=1$ puisque
$f(\varepsilon)\in({\bf Z}/b)^{\times(p)}$. on en d\'eduit
$\alpha\in ({\bf Z}/b)^{\times(p)}$, situation exclue.
Le sous-groupe $H(t)$ engendr\'e par  $t$ dans $Cl(A)$ est donc isomorphe \`a
${\bf Z}/n$. La fin de la d\'emonstration  est analogue \`a celle du 
th\'eor\`eme.
Les sous-groupes $H(t)$ et $H(t')$ engendr\'es respectivement par
$t=\partial([(\alpha+\sqrt{\delta})/2])$ et
$t'=\partial([(\alpha'+\sqrt{\delta})/2])$ sont en somme directe dans
$Cl(A)$.

Ces  \'el\'ements du groupe des classes ont \'et\'e construits pour
la premi\`ere fois par  Yamamoto ([28]).  \`A partir de ces
\'el\'ements, cet auteur  a montr\'e que pour tout
$n>1$, il existe une infinit\'e de corps quadratiques r\'eels et 
imaginaires dont le
groupe des classes contient un facteur ${\bf Z}/n$.
\endpreuve
\subsection{Construction d'\'el\'ements non triviaux de $Cl(A)_{(n)}$.}
Soit $F$ un corps de nombres d'anneaux d'entiers $A$.
Soit $r_1$ (resp. $2r_2$) le nombre de plongements
r\'eels (resp. complexes) de $F$.
Si $r=r_1+r_2-1$, on a
$rg(A^\times)=r$
et
$A^\times=K_1(A)=\mu\times\prod_{i=1}^r{\bf Z}\varepsilon_i$
o\`u $\mu$ est le groupe des racines de l'unit\'e contenues dans $A$ et
o\`u $\{\varepsilon_i, \ 1\leq i\leq r\}$
est un syst\`eme fondamental d'unit\'es de $A$.

Lorsque $A$ poss\`ede ``peu'' d'unit\'es,
l'extension $(\dagger)$ p. 20
permet d'obtenir des \'elements non triviaux de
$Cl(A)_{(n)}$ \`a partir d'\'el\'ements de
$K_1(A;{\bf Z}/n).$

L'\'egalit\'e $K_1(A; {\bf Z}/n)={\cal U}(A;{\bf Z}/n)$ et le lemme ($N$-$N_1$)
conduisent  au r\'esultat suivant:
\begin{prop}
Soit $F$ un corps de nombres, d'anneau d'entiers $A$ et soit $n$
un entier naturel. On suppose
qu'il existe $z\in A$ tel que $N(z)=b^n$,
$(N(z),N_1(z))=1$
et que pour
tout diviseur $m$ de $n$, $1\leq m<n$,
$b^m$ ne soit pas la norme d'un \'el\'ement
de
$F$.
Alors $[z]=z \ \mod\  F^{\times(n)}$
est un \'el\'ement d'ordre $n$ de
$K_1(A;{\bf Z}/n)$.

\noindent
Si $r_1+r_2-1=0$, on a $Cl(A)_{(n)}\not=0$.
Si $r_1+r_2-1\not=0$ et
si $-b$ n'est pas la norme
d'un \'el\'ement de $F$, alors
$Cl(A)_{(n)}\not= 0$.
\end{prop}
\preuve
d'apr\`es le lemme ($N$-$N_1$), $[z]$ appartient \`a
$K_1(A;{\bf Z}/n)$.
Si $[z]$ est d'ordre $m$,  $1\leq m<n$, il existe $u\in F^\times$
tel que
$z^m=u^n$, c'est-\`a-dire
$z=u^s$ avec $s=n/m$.
L'\'equation $N(z)=N(u)^s$
s'\'ecrit
$b^m=N(u)$, ce qui n'est pas. Notons comme toujours
$\partial:K_1(A;{\bf Z}/p)\to Cl(A)_{(n)}$ et
supposons \`a pr\'esent que
$\partial([z])=0$. Dans ce cas, il existe
$u\in F^\times$,
$\xi\in\mu$ et des entiers
$l_i$
tels que
$z=\xi^{l_0}\varepsilon_1^{l_1}\cdots\varepsilon_r^{l_r}u^n$,
d'o\`u l'on d\'eduit $N(z)=\pm N(u)^n$ (avec le signe $+$ si $r=0$)
soit
$b=\pm N(u)$, ce qui n'est pas.

\endpreuve

La classe caract\'eristique secondaire
$$d^{(n)}_1:Cl(A)_{(n)}\to \Omega^1_{dR}(A)/S$$
introduite en 1.6, corollaire 23 conduit au r\'esultat suivant.
\begin{prop}
Soit $F$ un corps de nombres d'anneaux d'entiers $A$.
On pose
$r=r_1+r_2-1$
et
$A^\times=\mu\times \prod _{i=1}^r{\bf Z}\varepsilon_i.$
Soit $n$ un diviseur du discriminant du corps $F$.

On suppose

\noindent 1) Pour tout
$\xi\in \mu$, $\xi^{-1}d\xi\equiv 0 \ \mod\  n\Omega^1_{dR}(A).$

\noindent 2) Pour tout
$i$, $1\leq i\leq r$, $\varepsilon_i^{-1}d\varepsilon_i\equiv 0 \ 
\mod\  n\Omega^1_{dR}(A).$

\noindent 3) Il existe $u\in A$ avec
$N(u)=b^n$, $(N(u),N_1(u))=1$, et

\noindent$u^{-1}du\not \equiv 0 \ \mod\  n\Omega^1_{dR}(A).$
Alors
$Cl(A)_{(n)}\not=0$.
\end{prop}
\preuve
les hypoth\`eses 1 et 2 montrent que
$dA^\times/A^\times\equiv 0\ \mod\  n\Omega^1_{dR}(A)$.
L'application
$d_{1}^{(n)}$,  de source
$Cl(A)_{(n)}$
est donc de but
$\Omega^1_{dR}(A)/(n)$.
D'apr\`es le lemme ($N$-$N_1$), les hypoth\`eses
3 fournissent l'\'el\'ement $[u]=u \ \mod\  F^{\times(n)}$ de
$K_1(A;{\bf Z}/n)$. De cet \'el\'ement,
on d\'eduit
$x=\partial([u])$ dans $Cl(A)_{(n)}$.
On a
$d_1^{(n)}(x)=u^{-1}du \  \mod \  n\Omega^1_{dR}(A)$,
quantit\'e non nulle par hypoth\`ese, ce qui montre que $x$ est
non trivial.
\endpreuve
\medskip
{\noindent
{\bf Exemple.} Posons $x={}^3\sqrt{182}$ et  soit $F={\bf Q}[x]$,
d'anneau d'entiers $A={\bf Z}[x]$,
d'unit\'e fondamentale $\varepsilon=17-3x$. Pour  $p=3$,
l'\'el\'ement $u=5-2x$ d\'efinit un \'el\'ement non nul de  $Cl(A)_{(p)}$.

  \subsection{Exemples de $n$-torsion du groupe des classes : cas d'un corps
quadratique imaginaire.}
\begin{prop}
Soit $F$ un corps quadratique d'anneau d'entiers $A$ et de 
discriminant $\delta<0$
et soit $n$ un entier impair.
On suppose qu'il existe $(\alpha,b)\in{\bf Z}^2$
tel que
$\alpha^2-4b^n=\delta$,
avec
$(\alpha, b)=1$.
On suppose de plus que pour tout diviseur
$m$ de $n$, $1\leq m<n$, la quantit\'e $\delta+4b^n$ n'est pas un 
carr\'e parfait.
Alors
$Cl(A)_{(n)}\not= 0.$
\end{prop}

\begin{rem}
Sous les m\^emes hypoth\`eses, si $\delta$ est positif et si de plus
$\pm b$ n'est pas une norme, la conclusion
$Cl(A)_{(n)}\not=0 $ subsiste.
\end{rem}
\preuve
on applique la proposition 37  \`a 
$z=\displaystyle\frac{\alpha+\sqrt\delta}{ 2}.$
L'\'equation $z^m=u^n$ conduit \`a
$\delta+4b^n$ carr\'e parfait.
\medskip

\medskip
Dans les quelques exemples ci-dessous, l'anneau $A$ des entiers du 
corps ${\bf Q}[\sqrt\delta]$
est tel que $Cl(A)_{(n)}\not=0$.

\noindent $n=3$

\noindent
$\delta=-104=2^2-4\cdot 3^3=4\cdot\ (-26)$

\noindent
$\delta=-5\ 320=2^2-4\cdot\ 11^3=4\cdot\ (-1330)$

\noindent
$\delta=-48\ 664=2^2-4\cdot\ 23^3=4\cdot\ (-12\ 166)$

\noindent $n=5$

\noindent
$\delta=-127=1^2-4\cdot\ 2^5$

\noindent
$\delta=-12\ 499=1^2-4\cdot\ 5^5$

\noindent
$\delta=-31\ 103=1^2-4\cdot\ 6^5$

\noindent
$\delta=-131\ 071=1^2-4\cdot\ 8^5$

\noindent
$\delta=-399\ 999=1^2-4\cdot\ 10^5$

\noindent $n=7$

\noindent
$\delta=-511=1^2-4\cdot\ 2^7$

\noindent
$\delta=-65\ 535=1^2-4\cdot\ 4^7$

\noindent
$\delta=-312\ 499=1^2-4\cdot\ 5^7$

\noindent $n=9$

\noindent
$\delta=-2047=1^2-2^9$

\noindent
$\delta=-78\ 728=2^2-3^9=4\cdot\ (-19\ 682)$

\noindent
$\delta=-78731=1^2-3^9$

\noindent $n=11$

\noindent
$\delta=-8191=1^2-2^{11}$

\noindent
$\delta=-708\ 584=2^2-3^{11}=4\cdot\ (-177\ 146)$

\noindent
$\delta=-708\ 587=1^2-3^{11}$

\subsection{Exemples de $n$-torsion ramifi\'ee du groupe des
classes: cas d'un corps quadratique}
Soit $F$ un corps
de nombres quadratique de discriminant $\delta$. Si $\delta<0$, on 
exclut les deux cas $\delta =-4$
et
$\delta=-3$ pour lesquels le groupe des classes est trivial
et le groupe des unit\'es n'est pas r\'eduit \`a ${\bf Z}/2$.

On pose $\omega=\frac{\sqrt{\delta}}{2}$ ou
$\omega=\frac{1+\sqrt{\delta}}{2}$
suivant que
$\delta\equiv\ 0\ \mod \ 4$ ou
$\delta \equiv 1\ \mod\ 4.$
L'anneau $A$ des entiers du corps $F$ est
${\bf Z}[\omega]$.
Posons $P=X^2-\delta $ si $\delta\equiv\ 0\ \mod \ 4$ et
$P=X^2-X+(1-\delta)/4$ sinon.
L'homologie de Hochschild de $A$ est donn\'ee par la
\begin{prop}

\item{a)}
Si $\delta\equiv 1\ \mod \ 4$, on a
$\Omega^1_{dR}(A)={\bf Z}/\delta \ d\omega$
et
$\omega d\omega=\frac{1}{2}d\omega$.
\item{b)}
Si $\delta\equiv 0\  \mod \ 4$, on a
$\Omega^1_{dR}(A)={\bf Z}/(\delta/2)\ d\omega \oplus{\bf Z}/2\  \omega
d\omega.$
\end{prop}

On en d\'eduit

\begin{prop}

\item{a)}
Soit $n$ un diviseur impair du discriminant $\delta$
du corps quadratique $F$.
Alors on a
$\Omega^1_{dR}(A)/(n)={\bf Z}/n d\omega$
avec $\omega d\omega=\frac{1}{2}d\omega$ si $\delta\ \equiv\ 1\ \mod\ 4$
et $\omega d\omega=0$ si $\delta\ \equiv \ 0\ \mod\ 4$.
\item{b)}
on suppose $\delta\ \equiv\ 0\ \mod\  4$ et $n$ diviseur pair de $\delta$.
Alors on a
$$\Omega^1_{dR}(A)/(n)={\bf Z}/n\  d\omega\oplus {\bf Z}/2\ \omega
d\omega.$$
\end{prop}

Tout \'el\'ement $z$ de $A$ s'\'ecrit $z=\displaystyle 
\frac{\alpha+\beta \sqrt\delta}{2}$
avec
$\alpha$ et $\beta$ dans $\bf Z$.
On note $N(z)$ sa norme,
$N_1(z)=tr(z)$ sa trace
et $\sigma(z)$ son conjugu\'e.
On a
$\sigma(z)=\displaystyle\frac{\alpha-\beta \sqrt\delta}{ 2}$,
$tr(z)=\alpha$ et
$N(z)=\displaystyle\frac{\alpha^2-\delta\beta^2}{4}$.
Dans
$\Omega^1_{dR}(A)$, on a la relation
$N(z)z^{-1}dz=\sigma(z)dz.$

Supposons que  $n$ soit un diviseur impair du discriminant $\delta$.
De
$N(z)\equiv \alpha^2/4 \ \ \mod\  n{\bf Z}$
et de
$\sigma(z)dz\equiv \frac{\alpha\beta}{2}\displaystyle d\omega \ \ \mod \
n\Omega^1_{dR}(A)$, on d\'eduit
que si
$(N(z),n)=1$, on a
$z^{-1}dz\equiv\displaystyle \frac{2\beta}{ \alpha}d\omega \ \ \mod \ 
n\Omega^1_{dR}(A).$

En particulier, si
$F$ est r\'eel  et si
$\varepsilon=\displaystyle\frac{\varepsilon_1+\varepsilon_2 \sqrt\delta}{ 2}$
est l'unit\'e fondamentale de $A$, on a toujours
$(\varepsilon_1,n)=1$ et donc
  $\varepsilon^{-1}d\varepsilon\equiv 0 \ \mod\ n\Omega^1_{dR}(A)$
si et seulement si
$n$ divise $\varepsilon_2$.
Dans ce cas, on a $dA^\times/A^\times \equiv\ 0\ \mod\ 
n\Omega^1_{dR}(A)$. On en d\'eduit
$$\Omega^1_{dR}(A)/(n,dA^\times/A^\times)=\Omega^1_{dR}(A)/(n)={\bf Z}/n\ dw.$$
On obtient la classe caract\'eristique
$$ d_1^{(n)}:Cl(A)_{(n)}\to {\bf Z}/n \ dw$$
Enfin, si $u=\displaystyle\frac{\alpha+\beta\sqrt\delta}{2}$
est un \'el\'ement de $A$ tel que
$(\beta, n)=(\alpha,n)=1$,
alors
$u^{-1}du\ \equiv \frac{2\beta}{\alpha} d\omega\ \mod\ n$
est une quantit\'e non nulle de ${\bf Z}/n\ d\omega$.
De tout ceci, on d\'eduit que la proposition 38 prend la forme:

\begin{prop}
Soit $F$ un corps quadratique de discriminant $\delta$ et
d'anneau d'entiers $A$. Soit $n$ un diviseur impair de
$\delta$. Si $F$ est r\'eel, on suppose que l'unit\'e fondamentale
$\varepsilon=\displaystyle\frac{\varepsilon_1+\varepsilon_2 \sqrt{\delta}}{ 2}$
est telle que $n|\varepsilon_2$.
Soient $(\alpha,\beta,b))\in{\bf Z}^3$ une solution de l'\'equation
$\alpha^2-4b^n=\delta\beta^2$ avec
$(b,\alpha)=(\beta, n)=(\alpha, n)=1.$
Alors
$Cl(A)$ poss\`ede un \'el\'ement d'ordre $n$.
\end{prop}

En se restreignant aux \'el\'ements $u$ de la forme
$\displaystyle\frac{\alpha+\sqrt{\delta}} {2}$, on obtient

\begin{prop}
Soient
$\alpha$, $b$ et $n$ trois entiers avec $n$ impair, $(\alpha, 
b)=(\alpha, n)=1$.
On pose $\delta=\alpha^2-4b^n$.
On suppose que
$n$ divise $\delta$ et que $\delta$ est le discriminant d'un corps
quadratique $F$ d'anneau d'entiers $A$.
Si $\delta$ est positif, on suppose de plus que
l'unit\'e fondamentale
$\varepsilon=
\displaystyle\frac{\varepsilon_1+\varepsilon_2\sqrt{\delta}}{2}$ de
$A$ est telle que $n|\varepsilon_2$.
Alors $Cl(A)_{(n)}\not= 0.$
\end{prop}

Dans les quelques exemples ci-dessous, l'anneau $A$ des entiers du
corps ${\bf Q}[\sqrt\delta]$
est tel que $Cl(A)$ poss\`ede un \'el\'ement d'ordre $n$.

\noindent $n=3$

\noindent
$\delta=231=17^2-4\cdot(-2)^3={\bf 3}\cdot 107$, d'unit\'e fondamentale
$\varepsilon=(430+24\sqrt{\delta})/2$.

\noindent
$\delta=-231=5^2-4\cdot\ 4^3=-{\bf 3}\cdot\ 7\cdot\ 11$

\noindent
$\delta=-255=1^2-4\cdot\ 4^3=-{\bf 3}\cdot\ 5\cdot\ 17$

\noindent
$\delta=-16\ 383=1^2-4\cdot\ 16^3=-{\bf  3}\cdot\ 43\cdot\ 127$

\noindent
$\delta=-62\ 484=4^2-4\cdot\ 25^3=4\cdot\ (-{\bf 3}\cdot\ 41\cdot\ 127)$

\noindent
$\delta=-3\ 999\ 999=1^2-4\cdot\ 100^3=-{\bf 3}\cdot\ 23\cdot\ 29\cdot\ 1999$

\noindent $n=5$

\noindent
$\delta=-236\ 195=1^2-4\cdot\ 9^5=-{\bf 5}\cdot\ 97\cdot\ 487$

\noindent
$\delta=-644\ 195=3^2-4\cdot\ 11^5=-{\bf 5}\cdot\ 19\cdot\ 6\ 781$

\noindent
$\delta=-9\ 904\ 380=4^2-4\cdot\ 19^2=4\cdot\ (-3\cdot\ {\bf 5}\cdot\ 
383\cdot\ 431)$

\noindent $n=7$

\noindent
$\delta=-511=1^2-4\cdot\ 2^7=-{\bf 7}\cdot\ 23$

\noindent
$\delta=-65\ 527=3^2-4\cdot\ 4^7=-{\bf 7}\cdot\ 11\cdot\ 23\cdot\ 37$

\noindent $n=11$

\noindent
$\delta=-708\ 587=1^2-4\cdot\ 3^{11}=-{\bf 11}\cdot\ 37\cdot\ 1741$

\noindent $n=15$

\noindent
$\delta=-4\ 294\ 967\ 295=1^2-4\cdot\ 4^{15}=-{\bf 3\cdot\ 5}\cdot\ 
17\cdot\ 257\cdot\ 65\ 537$

\medskip

\section{Applications \`a la cyclotomie.}
\subsection{ Notations et strat\'egie g\'en\'erale.}

Soit $p$ un nombre premier impair et soit $\zeta=\zeta_p$ une racine primitive
$p$-i\`eme de l'unit\'e. Le corps cyclotomique $F={\bf Q}[\zeta]$ est 
une extension
galoisienne de degr\'e $p-1$ de ${\bf Q}$, de groupe de Galois $G=\left({\bf
Z}/p\right)^\times$. Soit $g$ un g\'en\'erateur de $G$. On d\'esigne par $s$,
$1<s\leq p-1$, l'entier tel que
$g\zeta=\zeta^s$. La conjugaison complexe $g^{(p-1)/2}$ est not\'ee $\sigma$.
L'anneau $A$ des entiers de $F$ est
${\bf Z}[\zeta]$. Soient $Cl(A)$ le groupe des classes de $A$ et  $h=h_p$
le nombre de classes de $A$. Le sous-corps maximal r\'eel ${\bf 
Q}[\zeta+\zeta^{-1}]$ de
$F$ a pour nombre de  classes $h^+$. On sait que $h^+\mid h$ et que
$h^+$ est le nombre de classes de $A$ invariantes par conjugaison complexe.
La $p$-torsion du
groupe des classes $Cl(A)$ se d\'ecompose en
$Cl(A)_{(p)}=Cl(A)_{(p)}^-\oplus Cl(A)^+_{(p)}$
avec $Cl(A)^{\pm}_{(p)}=\ker(\sigma\mp\id)$.
On sait  que si  $p^a$ d\'esigne le nombre d'\'el\'ements de 
$Cl(A)^-_{(p)}$, alors
$p^a$ divise $h^-=h/h^+$.

La conjugaison complexe sur $A$ d\'efinit une involution toujours 
not\'ee $\sigma$ sur
$K_1(A;{\bf Z}/p)$ qui s'\'ecrit
$$K_1(A;{\bf Z}/p)=K_1^-(A;{\bf Z}/p)\oplus K_1^+(A;{\bf Z}/p)$$ avec
$K_1^\pm(A;{\bf Z}/p)=\ker (\sigma\mp\id)$. Par ailleurs,
$$A^\times/(p)\cong\ \mu_p\times\{\pm 1\}\times\left({\bf
Z}/p\right)^{(p-3)/2}$$
o\`u
$\mu_p=\{\exp(2ik\pi /p),\ 0\leq k\leq p-1\}$. En particulier,
$\left(A^{\times }/(p)\right)^-\cong\ \mu_p$. L'extension ($\dagger$) p. 20
se scinde donc en deux parties dont la partie antisym\'etrique s'\'ecrit

$$
\xymatrix{
1\ar[r] &\mu_p\ar[r] &K_1^-(A;{\bf Z}/p)\ar[r]^-{\partial} &Cl(A)^-_{(p)}\ar[r]
&1\cdot
  }
\leqno{(\dagger^-)}
  $$
Dans toute la suite de ce texte , on pose
$$d_p^-:=\dim_{{\bf Z}/p}Cl(A)_{(p)}^-=\dim_{{\bf Z}/p}K_1^-(A;{\bf Z}/p)-1.$$
Rappelons qu'un nombre premier est r\'egulier s'il ne divise pas le nombre de
classes $h_p$: pour un nombre premier r\'egulier,
$Cl(A)_{(p)}=Cl(A)^-_{(p)}=0$.

\begin{prop}
On a $d_p^-=0$ si et seulement si $p$ est un nombrer premier r\'egulier.
\end{prop}

\preuve
Si $p$ est r\'egulier, l'extension $(\dagger^-)$ se r\'eduit \`a
$K_1^-(A;{\bf Z}/p)=\mu_p$. R\'eci\-pro\-quement, si $K_1^-(A;{\bf
Z}/p)=\mu_p$, alors $h^-=0$ et $p$ ne divise pas $h^-$. D'apr\`es un 
th\'eor\`eme
de Kummer ([24], 5.6), ceci entra\^{\i}ne que $p$ ne divise pas $h^+$ donc
$Cl(A)_{(p)}=0$, c'est-\`a-dire $p$ r\'egulier.\endpreuve

\medskip
Rappelons que $(p,a,b,c)$ satisfont aux hypoth\`eses du premier cas du dernier
th\'eor\`eme de Fermat (en abr\'eg\'e DTF1) si $p$ est un nombre premier impair
et si
$a^p=b^p+c^p$ avec $(a,b,c)=(p,abc)=1$  (on parle du second cas si 
$p$ divise $abc$).

La d\'emarche d\'evelopp\'ee dans les paragraphes qui suivent est
celle-ci.
L'\'equation $a^p=b^p+c^p$ permet de construire un \'el\'ement $z$ de
$K_1^-(A;{\bf Z}/p)$. La trace de Dennis \`a coefficients nous permet
de montrer que cet \'el\'ement $z$ n'est pas trivial. L'action du 
groupe de Galois
fournit$(p-1)/2$ \'el\'ements de
$K_1^-(A;{\bf Z}/p)$ construits \`a partir de $z$. Gr\^ace \`a la 
trace de Dennis,
nous minorons la dimension du sous-espace vectoriel de $K_1^-(A;{\bf Z}/p)$
engendr\'e par ces $(p-1)/2$ \'el\'ements en termes des polyn\^omes de
Mirimanoff. On en d\'eduit une minoration de $d_p^-$.

\medskip
Au vocabulaire pr\`es, le r\'esultat suivant est bien connu.

\begin{prop}
Soient $p$ un nombre premier, $A$ l'anneau ${\bf Z}[\zeta_p]$, $F$ le corps
${\bf Q}[\zeta_p]$. On suppose que $(p,a,b,c)$ satisfont les hypoth\`eses du
premier cas du dernier th\'eor\`eme de Fermat.
Pour $1\leq \ell\leq (p-1)/2$, les \'el\'ements
$$z_\ell=\frac{a-b\zeta^{s^\ell}}{a-b\zeta^{-s^\ell}}\ \mod\ F^{\times(p)}$$
appartiennent alors \`a $K^-_1(A;{\bf Z}/p)$.
\end{prop}

\preuve Sous les hypoth\`eses DTF1, les id\'eaux fractionnaires principaux
$(a-b\zeta^\ell)$, $1\leq \ell\leq p-1$ sont deux \`a deux premiers 
entre eux. On
en d\'eduit que chacun de ces id\'eaux s'\'ecrit
sous la forme
$(a-b\zeta^\ell)=I^p_\ell$, o\`u les $I_\ell$ sont des id\'eaux 
fractionnaires. Par
cons\'equent,  pour $1\leq \ell\leq p-1$, les \'el\'ements 
$a-b\zeta^\ell\ \mod\
F^{\times(p)}$ appartiennent \`a
${\cal U}(A;{\bf Z}/p)$.
\endpreuve

\subsection{Emploi du groupe $K_1(A/p;{\bf Z}/p)$.}

Soit $\varphi:A\to A/p$ la projection canonique. Posons
$\lambda=1-\varphi(\zeta)$.
Alors $A/(p)={\bf Z}/p[\lambda]$ avec $\lambda^{p-1}=0$. L'anneau $A/p$ est
local et
$\left(A/p\right)^\times=\break{\bf Z}/p^\times+\lambda{\bf Z}/p[\lambda]$. On
en d\'eduit
$$K_1(A/p;{\bf Z}/p)=
\left(A/p\right)^\times/(p)=\left(1+\lambda{\bf Z}/p[\lambda],\times\right).$$

\medskip
Les modules de diff\'erentielles $\Omega^1_{dR}(A)/(p)$,
$\Omega^1_{dR}(A/p)$ et $\Omega^1_{dR}(A/p)/(p)$ sont tous trois isomorphes
\`a
$${\bf Z}/p[X]dX/(X-1)^{p-2}dX,$$
donc
$\Omega^1_{dR}(A/p)={\bf Z}/p[\lambda]d\lambda$
avec
$\lambda^{p-1}=0$ et $\lambda^{p-2}d\lambda=0$.

Par commodit\'e,  $o(\lambda^j)$ d\'esigne un \'el\'ement
ind\'etermin\'e de
$\lambda^{j+1}{\bf Z}/p[\lambda].$
Soit $p$ un nombre premier impair et soient $x$ et $y$ deux  \'el\'ements de
$({\bf Z}/p)^\times$ tels que $x-y=1$. Soint les \'el\'ements
$w=x-y(1-\lambda)$ et $\sigma(w)=x-y(1-\lambda)^{-1}$ de
$(A/p)^\times$ et soit
$z'=z'(x)$ l'\'el\'ement de $K_1^-(A/p;{\bf Z}/p)$ d\'efini par
$z'=w/\sigma(w)\ \mod\ (A/p)^{\times(p)}.$

\begin{prop}
Si $p$ est un nombre premier impair et si
\break $x\in{\bf Z}/p\setminus\{0,1,1/2\}$, alors l'\'el\'ement $z'(x)$
ci-dessus de
$K_1^-(A/p;{\bf Z}/p)$ n'est pas colin\'eaire \`a
l'\'el\'ement
$1-\lambda$.
\end{prop}

\preuve Calculons les traces $D^{(p)}_1(z'(x))$ et
$D_1^{(p)}(1-\lambda)$. On a
$$D_1^{(p)}(z'(x))= w^{-1}dw-\sigma(w)^{-1}d\sigma (w).$$
Puisque $w=1+y\lambda$,
$w^{-1}=\sum_{k\geq 0}(-1)^ky^k\lambda^k$,
$dw=yd\lambda$ et
$$w^{-1}dw=\sum_{k\geq 0}(-1)^ky^{k+1}\lambda^kd\lambda.$$
De
$\sigma(w)=\displaystyle\frac{1-\lambda x}{1-\lambda}$, on d\'eduit
$\sigma(w)^{-1}=1+\sum_{k\geq 1}x^{k-1}y\lambda^k$ tandis que
$$d\sigma(w)=-y\sum_{k\geq 1}k\lambda^{k-1}d\lambda$$ et par suite
$$\sigma(w)^{-1}d\sigma(w)
=-yd\lambda-y(y+2)\lambda
d\lambda-y(3+2y+xy)\lambda^2d\lambda+o(\lambda^2)d\lambda.$$ Ces
expressions de
$w^{-1}dw$ et $\sigma(w)^{-1}d\sigma(w)$ conduisent \`a
$$D_1^{(p)}(z'(x))=
2yd\lambda+
2y\lambda d\lambda+
(3y+3y^2+2y^3)\lambda^2d\lambda+o(\lambda^2)d\lambda.$$
Par ailleurs
$D_1^{(p)}(1-\lambda)=
-(1-\lambda)^{-1}d\lambda=
-\sum_{k\geq 0}\lambda^kd\lambda.$
Supposons $z'(x)$ et $1-\lambda$ colin\'eaires. La comparaison des
coefficients en $d\lambda$ et en $\lambda^2d\lambda$ des traces de Dennis
de $z'(x)$ et $1-\lambda$ conduit \`a  l'\'egalit\'e
$2y^3+3y^2+y=0.$
Puisque $y\not=0$, on en d\'eduit que $y\in({\bf Z}/p)^\times$ est solution de
l'\'equation
$2X^2+3X+1=0$ dans
${\bf Z}/p[X]$.
Ceci conduit \`a $y=-1$ ou $y=-1/2$. Or n\'ecessairement $y\not=-1$ car sinon
$x=0$, ce qui est exclu. Par ailleurs, $y=-1/2$ \'equivaut \`a 
$x=1/2$, situation
\'egalement exclue.\endpreuve

\begin{prop}
On suppose que $(p,a,b,c)$ satisfait DTF1 avec $p>3$.
Alors, $d_p^-\geq 1$.
\end{prop}
\preuve
D\'esignons par $\varphi_1:K_1(A;{\bf Z}/p)\to K_1(A/p;{\bf Z}/p)$
l'application induite par
$\varphi$ en $K$-th\'eorie \`a coefficients. L'\'el\'ement
$$z=\displaystyle\frac{a-b\zeta}{a-b\zeta^{-1}}\ \mod\ F^{\times(n)}$$
de $K_1^-(A;{\bf Z}/p)$ est tel que
$$\varphi_1(z)=\frac{x-y(1-\lambda)}{x-y(1-\lambda)^{-1}}\ \mod\
(A/p)^{\times(p)}$$ avec
$x=\overline a/\overline c$ et $y=x-1$. On a n\'ecessairement
$x\not=0$. Si
$x=1/2$,  l'\'el\'ement
$$z_1=\displaystyle\frac{a-c\zeta}{a-c\zeta^{-1}}\ \mod\ F^{\times(n)}$$
  est tel que $\varphi_1(z_1)=z'(x_1)$ avec $x_1=\overline c/\overline b$. Les
hypoth\`eses DTF1 montrent que pour $p>3$, il est impossible d'avoir
simultan\'ement
$x=x_1=1/2$. La proposition pr\'ec\'edente s'applique donc pour l'un des deux
\'el\'ements $z$ ou $z_1$.\endpreuve

On a remarqu\'e plus haut que $d_p^-=0$ caract\'erise les nombres premiers
r\'eguliers. On a donc montr\'e:

\begin{cor}(Kummer, 1847)

\noindent Soit $p$ un nombre premier r\'egulier. Alors le premier cas du
dernier th\'eor\`eme de Fermat est satisfait pour $p$.
\end{cor}

\subsection{Emploi du groupe $K_1(R;{\bf Z}/p)$ .}

  Dans tout ce paragraphe, $x$ et $y$  d\'esignent  deux
\'el\'ements de $({\bf Z}/p)^\times $ tels que $x-y=1. $

L'action du groupe de Galois
$G=Gal(F/{\bf Q})$ sur
$\Omega^1_{dR}(A)/(p)={\bf Z}/p[\lambda]d\lambda$ est peu lisible car
$g\lambda=1-(1-\lambda)^s$. C'est
pourquoi nous introduisons les anneaux
$R'={\bf Z}[X]/(X^p-1)={\bf Z}[t]$ et $R=R'/p$
avec $t=X\ \mod \ (X^p-1)$.
  Le groupe $G$ op\`ere sur $R'$ par $gt=t^s$ (o\`u $g$ est un g\'en\'erateur de
$G$ et $(s,p)=1$). Remarquons que l'involution $\sigma=g^{(p-1)/2}$ est
telle que
$\sigma(t)=t^{-1}$.

On a $R={\bf Z}/p[1-t]$ avec $(1-t)^p=0$. L'anneau $R$ est local,
$R^\times=({\bf Z}/p)^\times\oplus (1-t){\bf Z}/p[1-t]$ et
$$K_1(R;{\bf Z}/p)=(1+(1-t){\bf Z}/p[1-t],\times).$$

  Les modules de
diff\'erentielles de K\"ahler
$\Omega^1_{dR}(R')$,
$\Omega^1_{dR}(R')/(p)$, $\Omega^1_{dR}(R)$ et
$\Omega^1_{dR}(R)/(p)$ sont tous quatre isomorphes \`a
$${\bf Z}/p[X]dX/(X-1)^pdX.$$

L'action de $G$ sur $\Omega^1_{dR}(R)$ est donn\'ee par
$g(t^idt)=g(t)^idg(t)=st^{s(i+1)-1}dt.$
Pour $1\leq k\leq p-1$, les relations
$$g(t^{s^k}t^{-1}dt)=st^{s^{k+1}}t^{-1}dt,\ \
\sigma(t^{s^k}t^{-1}dt)=-t^{-s^k}t^{-1}dt$$
et
$$g(t^{-1}dt)=st^{-1}dt,\ \ \sigma(t^{-1}dt)=-t^{-1}dt$$
conduisent \`a la d\'ecomposition commode suivante.
\begin{prop}
Posons $f^-_0=t^{-1}dt$,
et
pour $1\leq \ell\leq (p-1)/2$,\break
$f^\pm_\ell=\left(t^{s^\ell}\mp t^{-s^\ell}\right)t^{-1}dt.$

  On a alors
  $$\Omega^1_{dR}(R)=\Omega^-_{dR}(R)\oplus\Omega^+_{dR}(R)$$
  o\`u $\Omega^-_{dR}(R)$ est de dimension
  $(p+1)/2$, de base $(f^-_0,f^-_1,\cdots,f^-_{(p-1)/2})$
  et o\`u
   $\Omega^+_{dR}(R)$ est de dimension $(p-1)/2$, de base
   $(f^+_1, \cdots,f^+_{(p-1)/2})$.

   De plus, en d\'esignant par $g$ un g\'en\'erateur du groupe de Galois
   $G=Gal(F/{\bf Q})$ et en notant $\sigma$ l'involution
   $g^{(p-1)/2}$, on a les relations
   $g(f_0^-)=sf^-_0$,\ \  $g(f^\pm_\ell)=sf^\pm_{\ell+1}$,
   $1\leq \ell<(p-1)/2$,\ \ $g(f^{\pm}_{(p-1)/2})=f_1^\pm$

\noindent et
$\sigma(f^-_0)=-f^-_0$,\ \ $\sigma(f^\pm_\ell)=\pm f^\pm_\ell$
$(1\leq \ell\leq (p-1)/2$.
\end{prop}

\begin{defi}
Pour $x\in{\bf Z}/p\setminus\{0,1\}$, on pose $y=x-1$ et
  pour $1\leq k\leq (p-1)/2$, on introduit les \'el\'ements
$$\alpha_k=(x/y)^{s^{k-1}}+(y/x)^{s^{k-1}}$$ de
${\bf Z}/p$ et les \'el\'ements suivants de $K_1(R;{\bf Z}/p)$:
$v_k(x)=x-yt^{s^k}\ \mod\  R^{\times(p)}$,
$\sigma(v_k(x))=x-yt^{-s^k}\ \mod\  R^{\times(p)}$ et
$$z_k(x)=v_k(x)/\sigma(v_k(x)).$$
\end{defi}

\begin{prop}
Dans la base $(f_0^-,\cdots, f_{(p-1)/2}^-)$ de $\Omega^-_{dR}(R)$, la trace
de \break Dennis de $z_1(x)$ s'\'ecrit
$$D_1^{(p)}(z_1(x))=-s(x-1)\left( 2f_0^-+\sum_{k=1}^{(p-1)/2}\alpha_kf_k^-
\right).$$
\end{prop}
\preuve
On a $D_1^{(p)}(z_1(x))=
v_1^{-1}(x)dv_1(x)-\sigma(v_1(x))^{-1}d\sigma(v_1(x))$. Pour calculer
$v_1^{-1}(x)$, \'ecrivons $v_1(x)=-yt^s(1-(x/y)t^{-s})$.
L'identit\'e
$$\left(1-(x/y)t^{-s}\right)
\left(
1+(x/y)t^{-s}+\cdots+(x/y)^{p-1}t^{-(p-1)s}
\right)=1-x/y=-1/y$$ conduit
\`a
$$v_1^{-1}(x)=t^{-s}\left(1+(x/y)t^{-s}+\cdots+(x/y)^{p-1}t^{-(p-1)s}\right).$$
Puisque
$dv_1(x)=-syt^{s}t^{-1}dt$, on obtient
$$v_1^{-1}(x)dv_1(x)=
-sy\left( t^{-1}dt+\sum_{i=1}^{p-1}(x/y)^it^{-is}t^{-1}dt\right).$$
On transforme cette quantit\'e en \'ecrivant
$$v_1^{-1}(x)dv_1(x)=-sy\left(
t^{-1}dt
+\sum_{k=1}^{p-1}(x/y)^{s^{k-1}}t^{s^k}t^{-1}dt
\right)$$
soit encore
$$v_1^{-1}(x)dv_1(x)=-sy\left(
t^{-1}dt+
\sum_{k=1}^{(p-1)/2}(x/y)^{s^{k-1}}t^{s^k}t^{-1}dt
+\sum_{k=1}^{(p-1)/2}(x/y)^{-s^{k-1}}t^{-s^k}t^{-1}dt
\right).$$
Pour obtenir l'expression de $v_1^{-1}(x)dv_1(x)$ dans la base propos\'ee de
$\Omega^-_{dR}(R)$, introduisons
$\beta_k=(x/y)^{s^{k-1}}-(y/x)^{s^{k-1}}$. On a
$$v_1^{-1}(x)dv_1(x)=
-sy\left(
f_0^-+\frac{1}{2}\sum_{k=1}^{(p-1)/2}\alpha_k f^-_k
+\frac{1}{2}\sum_{k=1}^{(p-1)/2}\beta_k f^+_k
\right).$$
Le calcul de $\sigma(v_1(x))^{-1}d\sigma(v_1(x))$ se d\'eduit
imm\'ediatement de cette derni\`ere relation car
$D_1^{(p)}$ est \'equivariante,
$\sigma(f^-_0)=-f^-_0$, $\sigma(f^\pm_k)=\pm f^\pm_k$.
On obtient ainsi
$$\sigma(v^{1}(x))^{-1}d\sigma(v_1(x))=-sy\left(
-f_0^--\frac{1}{2}\sum_{k=1}^{(p-1)/2}\alpha_k f^-_k
+\frac{1}{2}\sum_{k=1}^{(p-1)/2}\beta_k f^+_k
\right),$$
d'o\`u finalement
$$D_1^{(p)}(z_1(x))=
-sy
\left(
2f_0^-+\sum_{k=1}^{(p-1)/2}\alpha_k f^-_k
\right).$$
\endpreuve

\begin{defi} On note $V(x)$ le sous-espace vectoriel  de $K_1^-(R;{\bf Z})$
engendr\'e par l'orbite de $z_1(x)$ sous l'action du groupe de Galois $G$,
c'est-\`a-dire
$$V(x)={\hbox{Vect}_{{\bf Z}/p}}(z_k(x),\ 1\leq k\leq (p-1)/2).$$
\end{defi}

\begin{prop}
Soit $C=C(x)$ la matrice circulante d'ordre $(p-1)/2$ \`a coefficients
dans ${\bf Z}/p$

$$C=C(x)=
\left(
\begin{array}{cccc}
\alpha_1, & \alpha_2, &\ldots , &  \alpha_{\frac{p-1}{2}} \\
  \alpha_{\frac{p-1}{2}}, & \alpha_1, & \ldots ,&  \alpha_{\frac{p-3}{2}}\\
\vdots & \ddots & \ddots & \vdots \\
\alpha_2, & \alpha_3, &\ldots,  & \alpha_1
\end{array}
\right)
$$
Alors
$$\dim_{{\bf Z}/p}V(x)\geq \hbox{rg}(C(x)).$$
\end{prop}
\preuve
\`a une constante pr\`es, les composantes de $D_1^{(p)}(z_1(x))$ dans la base
$(f_0^-,f^-_1,\cdots,f^-_{(p-1)/2})$ de
$\Omega^-_{dR}(R)$ sont
$(2,\alpha_1,\cdots,\alpha_{(p-1)/2}).$
Puisque $z_k(x)=g^k(z_1(x))$, compte tenu de l'action de $g$ sur les vecteurs
de base\break $(f_0^-,f^-_1,\cdots,f^-_{(p-1)/2})$, on en d\'eduit 
que la matrice
des composantes respectives de
$D_1^{(p)}(z_1(x))$, $D_1^{(p)}(z_2(x))$, $\cdots$, $D_1^{(p)}(z_{(p-1)/2}(x))$
a le m\^eme rang que la matrice
$$\left(
\begin{array}{ccccc}
2&\alpha_1, & \alpha_2, &\ldots , &  \alpha_{\frac{p-1}{2}} \\
  2&\alpha_{\frac{p-1}{2}}, & \alpha_1, & \ldots ,&  \alpha_{\frac{p-3}{2}}\\
2&\vdots & \ddots & \ddots & \vdots \\
2&\alpha_2, & \alpha_3, &\ldots,  & \alpha_1
\end{array}
\right).
$$
Puisque $\sum_{k=1}^{(p-1)/2}\alpha_k=-1$, le rang de cette matrice
est celui de la matrice $C(x)$. L'image de
$V(x)$ par la trace de Dennis $D_1^{(p)}$ a pour dimension le rang de
$C(x)$, cqfd.

\medskip
Le calcul du rang de la matrice $C(x)$ n\'ecessite l'introduction des
polyn\^omes de Mirimanoff.

\begin{defi}
Les polyn\^omes de Mirimanoff $M_k(X)\in{\bf Z}/p[X]$ sont d\'efinis pour
$1\leq k\leq p-1$ par
$$M_k(X)=\sum_{j=1}^{p-1}j^{k-1}X^j.$$
\end{defi}
Pour $t\in{\bf Z}/p$, on pose
$$r_p(t)=\#\{k\mid 1\leq k\leq (p-1)/2,\ M_{2k+1}(t)\not=0\}\cdot$$
C'est le nombre de  polyn\^omes de
Mirimanoff $M_j(X)$  non nuls en la valeur $t$ et d'indice $j$ impair.

\begin{prop} Soient $x$ et $y$ deux \'el\'ements de $({\bf 
Z}/p)^\times $  tels que
$x-y=1$.
Les valeurs propres de la matrice $C(x)$ sont $M_{2k+1}(x/y)$, $1\leq k\leq
(p-1)/2$. Le rang de la matrice $C(x)$  est
$r_p(x/y)$.
\end{prop}

\preuve
Soit $s$ le g\'en\'erateur de
$\left({\bf Z}/p\right)^\times$
qui d\'etermine l'action du groupe de Galois $G$ sur $A$ et soit
$v=s^2$ le g\'en\'erateur de
${\bf Z}/(p-1)/2\subset \left({\bf Z}/p\right)^\times$.
Les valeurs propres de la matrice $C$ sont alors

$$\aligned
\lambda_k &=\sum_{j=1}^{(p-1)/2}\alpha_j(v^k)^{j-1}\\
&=\sum_{j=1}^{(p-1)/2}(x/y)^{s^{j-1}}\left(s^{j-1}\right)^{2k}
+
\sum_{j=1}^{(p-1)/2}(x/y)^{s^{j-1+(p-1)/2}}\left(s^{j-1+(p-1)/2}\right)^{2k}
\\
&=\sum_{j=1}^{p-1}j^{2k}(x/y)^j\\
&=M_{2k+1}(x/y)\\
\endaligned
$$
Le rang de $C(x)$ est le nombre de valeurs propres non nulles.
Ces valeurs propres \'etant les $M_{2k+1}(x/y)$, le rang de
$C(x)$ est bien $r_p(x/y)$.
\endpreuve
En r\'esum\'e, nous avons montr\'e:

\begin{theo}
Soient $x$ et $y$ deux \'el\'ements de $({\bf Z}/p)^\times $  tels que
$x-y=1$. Alors
$$\dim_{{\bf Z}/p}V(x)\geq r_p(x/y).$$
\end{theo}
\rem
Posons
$$r_p=min\{r_p(t),\ t\in{\bf Z}/p\setminus\{0,1,1/2\}\}.$$
Alors, pour tout $x\in{\bf Z}/p\setminus\{0,1,1/2\}$, on a
$$(p-1)/2\geq \dim_{{\bf Z}/p}V(x)\geq r_p.$$
\subsection{Lien avec les d\'eriv\'ees logarithmiques de Kummer.}

Soit toujours  $A$ l'anneau des entiers du corps cyclotomique $F={\bf
Q}[\zeta_p]$ avec $p$ premier impair. On pose $\lambda=1-\zeta$. Identifions
$K_1(A/p;{\bf Z}/p)$ au groupe multiplicatif
$(1+\lambda{\bf Z}/p[\lambda],\times).$
Dans ses recherches sur le dernier th\'eor\`eme de Fermat pour les nombres
premiers irr\'eguliers, Kummer a introduit certaines ``d\'eriv\'ees
logarithmiques''.  Un \'el\'ement $z=\sum_{i=0}^{p-1}a_i\zeta^i$ de $A$, non
divisible par
$1-\zeta$ d\'etermine un \'el\'ement de $K_1(A/p;{\bf Z}/p)$ encore not\'e $z$.
Pour $1\leq k\leq p-2$, la d\'eriv\'ee logarithmique $\ell_k(z)$ est d\'efinie
comme la classe modulo $p$ de l'entier
$$\frac{d^k}{dX^k}\left(
\hbox{log}\left(\sum_{i=0}^{p-2}a_ie^{iX}\right)\right)_{X=0}.$$
Kummer a montr\'e que
$$\ell_k:(K_1(A/p;{\bf Z}/p),\times)\to ({\bf Z}/p,+)$$
est un morphisme de groupes.

Soient $x$ et $y$ deux \'el\'ements de $({\bf Z}/p)^\times$ tels que
$x-y=1$. L'\'el\'ement $z'(x)=\displaystyle\frac{x-y\zeta}{x-y\zeta^{-1}}$ de
$K_1^-(A/p;{\bf Z}/p)$ est tel que
$\ell_{2k}(z'(x))=0$ , $\ell_{2k+1}(z'(x))=2\ell_{2k+1}(x-y\zeta).$

Mirimanoff a montr\'e ({\sl cf}. [20], VII ou [9]) que pour
$1\leq k\leq (p-3)/2$, on a l'\'egalit\'e
$\ell_{2k+1}(x-y\zeta)=-xM_{2k+1}(x/y)$.
Ceci permet de formuler un lien entre la trace de Dennis
\`a coefficients et les d\'eriv\'ees logarithmiques de Kummer. Soient
$$z_k(x)=\frac{x-yt^{s^k}}{x-yt^{s^{-k}}}\  \mod\ (R)^{\times (p)}$$
les \'el\'ements de $K_1(R;{\bf Z}/p)$ introduits \`a la d\'efinition 51. Soit
  $C(x)\in \hbox{Mat}_{(p-1)/2}({\bf Z}/p)$ la matrice des coordonn\'ees de
$D_1^{(p)}(z_1(x))$, $D_1^{(p)}(z_2(x))$, ..., $D_1^{(p)}(z_{(p-1)/2}(x))$
dans la base de $\Omega^-_{dR}(R)$ d\'ecrite dans la proposition 50. Alors,
\`a une cons\-tante pr\`es, la matrice $C(x)$ a pour valeurs propres les
d\'eriv\'ees loga\-rith\-miques de Kummer
$\ell_{2k+1}(x-y\zeta)$.

Signalons un autre lien entre les $\ell_{2k+1}(x-y\zeta)$ et la trace 
de Dennis. Le
d\'eveloppement limit\'e
\`a l'ordre 2 de $D_1^{(p)}(z'(x))$ effectu\'e \`a la proposition 47 
peut \^etre
pr\'ecis\'e. On obtient
$$D_1^{(p)}(z'(x))=\sum_{k=0}^{p-3}\gamma_k(x)\lambda^kd\lambda$$
avec $\gamma_0(x)=2y$ et
$$\gamma_k(x)=(-1)^ky^{k+1}+(k+1)y+\sum_{j=1}^kjy^2(1+y)^{k-j}.$$
Introduisons les vecteurs colonnes  $\ell(x)$ et $D(x)$ de $({\bf
Z}/p)^{(p-1)/2}$ d\'efinis par
$\ell(x)=\left(\ell_{2k+1}(z'(x))\right)_{1\leq k\leq (p-1)/2}$ et
$D(x)=\left( \gamma_{2k}(x)\right)_{0\leq k\leq (p-3)/2}$.
Pour $p\leq 13$, on {\sl constate} qu'il existe une matrice triangulaire
$A\in \hbox{\bf GL}_{(p-1)/2}({\bf Z}/p)$ telle que
$\ell(x)=AD(x)$ pour tout $x\in ({\bf Z}/p)^\times$, ce qui montre qu'il est
\'equivalent de conna\^{\i}tre la trace de Dennis \`a coefficients ou les
d\'eriv\'ees logarithmiques de Kummer. Il serait int\'eressant de savoir si
cette observation se g\'en\'eralise \`a tout nombre premier.

\subsection{Application au premier cas du dernier th\'eor\`eme de Fermat.}

Supposons que $(p,a,b,c)$ satisfont aux hypoth\`eses DTF1.  Notons
$\overline a$, $\overline b$ et $\overline c$ les classes respectives 
de $a$, $b$
et $c$ dans ${\bf Z}/p$. Introduisons le sous-espace vectoriel
$V(p,a,b,c)$ de $K_1^-(A;{\bf Z}/p)$ engendr\'e par l'orbite de
$$z=z_1=\frac{a-b\zeta^s}{a-b\zeta^{-s}}\ \mod\ F^{\times(p)},$$
c'est-\`a-dire
$$V(p,a,b,c)=\hbox{Vect}_{{\bf Z}/p}(z_k,\ 1\leq k\leq (p-1)/2).$$

\begin{prop}On pose $x=\overline a/\overline c$ et
$y=1-x=\overline b/\overline
c$.
  Avec les notations de la section pr\'ec\'edente, on a
$$1\geq \dim_{{\bf Z}/p}V(x)-\dim_{{\bf Z}/p}V(p,a,b,c)\geq 0.$$
\end{prop}

\preuve soient $\varphi:A\to A/p$ la surjection canonique et  $\psi:R\to
A/p$ le morphisme d'anneaux d\'efini par  $\psi(t)=1-\lambda$. On d\'esigne
par $\varphi_1:K_1(A;{\bf Z}/p)$ et $\psi_1:K_1(R;{\bf Z}/p)\to
K_1(A/p;{\bf Z}/p)$ les applications induites en $K$-th\'eorie \`a
coefficients. Dans
$K_1(A;{\bf Z}/p)$, l'image de $V(p,a,b,c)$ par $\varphi_1$ co\"{\i}ncide
avec l'image de $V(x)$ par $\psi_1$, ce qui montre que la dimension de
$V(p,a,b,c)$ est sup\'erieure \`a celle de $\psi_1(V(x))$. On v\'erifie
ais\'ement que $\psi_1$ est surjective de noyau de dimension $1$. On en
d\'eduit l'in\'egalit\'e propos\'ee.
\endpreuve
\begin{theo}
Soient $(p,a,b,c)$ des entiers satisfaisant aux hypoth\`eses DTF1. Alors,
on a les in\'egalit\'es
$$d_p^-\geq r_p(\overline a/\overline c)-2\geq r_p-2.$$
\end{theo}
\preuve D'apr\`es le th\'eor\`eme 50 et la proposition ci-dessus, on a les
in\'egalit\'es
$$\displaylines{d_p^-\geq \dim_{{\bf Z}/p}K_1^-(A;{\bf Z}/p)-1\geq
\dim_{{\bf Z}/p}V(p,a,b,c)-1\geq\hfill\cr
\hfill\dim_{{\bf Z}/p}V(\overline a/\overline c)-2\geq
r_p(\overline a/\overline c)-2\geq r_p-2.\cr}$$

\rem
A normalisation pr\`es, les calculs ci-dessus correspondent \`a ceux 
effectu\'es
par Br\" uckner ([6]). Soit $\chi'$ la restriction de la trace de 
Dennis $D_1^{(p)}$
\`a l'espace $V(\overline a/\overline c)$. Notre trace $\chi'$ est \`a comparer
avec le morphisme
$\chi$ de [6],2.1. Les quantit\'es  $f_i(\eta)$ introduites en [6], 
3.5 sont telles
que
$f_i(\eta)\cong(-1)^{i-1}yM_{i-1}(\overline a/\overline c)\ \mod\  p$ et la
minoration
$d_p\geq r_p-2$ cor\-res\-pond \`a l'in\'egalit\'e [6], 5.1. \`A 
partir de cette
minoration, Br\" uckner montre que le premier cas du dernier th\'eor\`eme de
Fermat est vrai si
$p\geq 2^{d_p+3}-2d_p-3$, o\`u $ d_p=\dim_{{\bf Z}/p}Cl(A)_{(p)}.$ On 
peut aussi
exploiter l'in\'egalit\'e $d_p^-\geq r_p-2$ en proc\'edant comme suit.

\begin{prop}
Soit $p$ un nombre premier.
  On a $$d_p^-<\frac{p+3}{4}$$
\end{prop}

\preuve
LA quantit\'e  $p^{d_p^-}$
divise $h^-$. D'apr\`es [16] et [18], on a
$$h^-\leq 2p\left({\frac{p}{24}}\right)^{\frac{p-1}{4}}.$$
On en d\'eduit $$d_p-\frac{p+3}{4}\leq
\frac{\ln(2)}{\ln(p)}-\frac{(p-1)\ln(24)}{4\ln(p)}.$$
Le second membre de cette in\'egalit\'e est n\'egatif pour $p\geq 2$.

\medskip
De l'in\'egalit\'e
$d_p^-<(p+3)/4$ valable pour tout $p$ et de l'in\'egalit\'e
$d_p^-\geq r_p-2$, conditionnelle \`a une solution \`a DTF1, on d\'eduit le
r\'esultat suivant.

\smallskip
\noindent
{\bf Scholie}
{\sl Soit $p\geq 3$ un nombre premier. Si $r_p\geq (p+11)/4$, alors le premier
cas du dernier th\'eor\`eme de Fermat est satisfait pour $p$.}

\smallskip
Soulignons que le calcul num\'erique  de $r_p$ est assez rapide, ce 
qui ne semble
pas \^etre le cas pour $d_p$ ou $d_p^-$. Pour
$p<1000$, un calcul sur ordinateur montre qu'on a toujours  l'in\'egalit\'e
$r_p>(p+11)/4$ (le nombre maximal de valeurs nulles pour
$M_{2k+1}(t)$ est
$7$).

\subsection{Lien avec les nombres de Bernoulli.}
L'in\'egalit\'e $d_p^-\geq r_p-2$ propos\'ee
au th\'eor\`eme 60  peut se retrouver par un autre raisonnement.
  Le nombre
$r_p$ est reli\'e \`a la divisibilit\'e des nombres de Bernoulli au 
moyen des congruences
de Kummer. Rappelons en premier lieu  que les nombres de Bernoulli $B_k\in{\bf
Q}$ sont d\'efinis par
$$\frac{X}{\exp(X)-1}=\sum_{k\geq 0}B_k\frac{X^k}{k!}.$$
Soit $i(p)$ l'indice d'irr\'egularit\'e de $p$ d\'efini comme le nombre de
nombres de Bernoulli divisibles par $p$ (c'est-\`a-dire dont le num\'erateur
est divisible par $p$).  On a $$i(p)=\# \{ k,\ 1\leq k\leq (p-3)/2,\ p\mid
B_{2k}\}.$$

Rappelons en second lieu que pour  $x\in{\bf Z}/p\setminus\{0,1\}$, on dit que
$x$ satisfait les congruences de Kummer $({\cal K})$ si
$$B_{p-(2k+1)}M_{2k+1}(x)\cong 0\ \ \ \mod\ p\ \ (1\leq k\leq
(p-3)/2).\leqno{({\cal K})}$$

Il est clair que si $x$ est solution des congruences de Kummer, on a
l'in\'egalit\'e $$r_p(x)\leq i(p).$$
Par ailleurs, Kummer a montr\'e
que si $(p,a,b,c)$ satisfont aux hypoth\`eses DTF1, alors $x=\overline
a/\overline c$ satisfait les congruences $({\cal K})$ ({\sl cf}. 
[20], VII ou [9]). On
en d\'eduit $r_p\leq i(p)$. Cette in\'egalit\'e est conditionnelle 
\`a l'existence
d'une solution \`a DTF1. L'in\'egalit\'e  $i(p)\leq d_p^-$, ind\'ependante
d'une \'eventuelle solution \`a DTF1, r\'esulte d'un th\'eor\`eme de 
Ribet [21].

\bigskip

\noindent {\bf Bibliographie}
\bigskip
\medskip

\noindent
[1], Artin E. \& Tate J., {\sl Class field theory}, Benjamin, New-York, 1967.

\medskip
\noindent
[2] Bass H., {\it{Algebraic $K$-theory}}, Benjamin,
New-York, 1968.

\medskip
\noindent
[3] Bass H., Milnor J.  \& Serre J.-P., Solution of the congruence subgroup
pro\-blem for $SL_n$ ($n\geq 3$) and $Sp_{2n}$ ($n\geq 2$), {\sl Publ. Math.
Inst. Hautes \'Etudes Sci.}, {\bf 33}, 1967, pp. 59-37.
\medskip

\noindent
[4] Berrick J., Interwiners and the $K$-theory of commutative rings,
{pr\'epublication}, 2000.

\medskip
\noindent
[5] Bourbaki N., {\it{Alg\`ebre}}, chap. 1-3, Hermann,
Paris, 1970.

\medskip
\noindent
[6], Br\" uckner H., Zum ersten Fall der Fermatschen Vermutung,
{\sl J. Reine ang. Math.}, {\bf 274-276}, 1975, pp. 21-26.

\medskip
\noindent
[7] Connes A., Non-commutative differential geometry,
{\it Publ. Math. Inst. Hautes \'Etudes Sci.}, {\bf 62}, 1985, pp. 257-360.

\medskip
\noindent
[8]  Goodwillie T.-G., Relative algebraic $K$-theory and cyclic homology, {\sl
Ann. of Math.}, {\bf 124}, 1986, pp. 347-402.

\medskip
\noindent
[9] Granville A., The Kummer-Wieferich-Skula approach to the first case of
Fermat's last theorem, Proceedings of the $3^{rd}$ conference of the
Canadian Number theory Association, August 18-25 1991, {\sl Advances in
Number theory}, ed. F.-Q. Gouve\^a \& N. Yui, Clarendon Press, Oxford, 1993.

\medskip
\noindent
[10]  Igusa K., What happens to Hatcher and Wagoner's formula for $\pi_0C/M$
when the first Postnikov invariant of $M$ is trivial?, {\sl Lectures Notes in
Math.}, {\bf 1046}, New-York, Springer Verlag, 1984,  pp. 104-72.

\medskip
\noindent
[11]  Karoubi M., Homologie cyclique et $K$-th\'eorie,
{\it{Ast\'erisque}}, {\bf 149}, Soc. Math. France, 1987.
\medskip

\noindent
[12],  Karoubi M. \&  Lambre T., Quelques classes caract\'eristiques 
en th\'eorie
des nombres, {\sl C. R. Acad. Sci. Paris}, t. 330, S\'erie I, 2000.

\medskip
\noindent
[13]  Karoubi M. \&  Villamayor O., $K$-th\'eorie alg\'ebrique et $K$-th\'eorie
topo\-lo\-gique I, {\sl Math. Scand.}, {\bf 28}, 1971, pp. 265-307.

\medskip
\noindent
[14]  Lambre T., Quelques exemples de lemme de premi\`ere
perturbation en homologie cyclique, {\it{Comm. Algebra}}, {\bf
23}, 1995, pp. 525-541.

\medskip
\noindent
[15]  Larsen M.,  Lindenstrauss A., Cyclic homology of Dedekind domains,
{\it{$K$-theory}}, {\bf 6}, 1992, pp. 301-334.

\medskip
\noindent
[16]  Lepist\"o T., On the growth of the first factor of the class
number of the prime cyclotomic field, {\it{Ann. Acad. Sci.
Fennicae}}, S\'erie~A, I, {\bf {577}}, 1974, Helsinski (21 pages).

\medskip
\noindent
[17]  Loday J.-L., {\sl Cyclic Homology}, Springer Verlag, Berlin, 1992.

\medskip
\noindent
[18]  Mets\"ankyl\"a T., Class numbers and $\mu$-invariants of
cyclotomic fields, {\it{Proc. Amer. Math. Soc.}}, {\bf 43}, 2, 1974,
pp. 299-300.

\medskip
\noindent
[19]  Neisendorfer J., Primary homotopy theory,
{\sl Memoirs Am. Math. Soc.}, {\bf 232}, 1980.

\medskip
\noindent
[20] Ribenboim P., {\it{13 lectures on Fermat's Last
Theorem}}, Springer, Berlin, 1974.

\medskip
\noindent
[21],  Ribet K., A modular construction of unamified $p$-extensions of ${\bf
Q}(\mu_p)$,
{\sl Invent. Math.}, {\bf 34}, 1976, pp. 151-162.

\medskip
\noindent
[22]  Serre J.-P., {\it{ Corps locaux}}, Hermann, 1968.

\medskip
\noindent
[23]  Wagoner J.-B., Delooping classifying spaces in algebraic $K$-theory,
{\sl Topology}, {\bf 1972}, 11, pp. 349-370.

\medskip
\noindent
[24]  Washington L., {\it{Introduction to cyclotomic
Fields}}, GTM 83, Springer, Berlin, 1982 .

\medskip
\noindent
[25]  Weibel Ch., Nil $K$-theory maps to cyclic homology, {\sl Trans. Am.
Math. Soc.}, {\bf 303}, 1987, pp. 541-558.

\medskip
\noindent
[26]  Weibel Ch., {\sl An introduction to homological algebra}, Cambridge
Studies in advanced mathematics, {\bf 38}, Cambrigde, 1994.
\medskip

\noindent
[27]  Wodzicki M., Excision in cyclic homology and in rational algebraic
$K$-theory, {\sl Ann. Math.}, {\bf 129}, 1989, pp. 591-639.

\medskip
\noindent
[28]  Yamamoto Y., On unramified Galois extensions of quadratic
number fields, {\it{Osaka J. Math.}}, {\bf 7}, 1970, pp. 57-76.

\end{document}